\documentclass[preprint]{elsarticle}

\usepackage[body={6in,9in},top=1in,left=1in,dvips]{geometry}

\biboptions{sort&compress}

\usepackage{amsfonts,amsmath,graphicx,subfigure,caption,bm}
\usepackage{xcolor}
\usepackage{algorithm,algorithmicx,algpseudocode}

\usepackage{float}
\usepackage{url}

\newcommand{\M}{\mathcal{M}}

\newcommand{\Sph}{\mathbb{S}}

\newcommand{\ds}{\displaystyle}

\newcommand{\vect}[1]{\mathbf{#1}}

\newcommand{\vx}{\vect{x}}
\newcommand{\xh}{\hat{x}}
\newcommand{\yh}{\hat{y}}
\newcommand{\vxh}{\vect{\xh}}
\newcommand{\vy}{\vect{y}}

\newcommand{\vf}{\vect{f}}

\newcommand{\vu}{\vect{u}}

\newcommand{\vxi}{\boldsymbol{\xi}}
\newcommand{\veta}{\boldsymbol{\eta}}

\newcommand{\X}{\mathbf{X}}
\newcommand{\Xh}{\hat{\X}}
\newcommand{\uc}{\underline{c}}

\newcommand{\ub}{\underline{b}}
\newcommand{\ua}{\underline{a}}
\newcommand{\uu}{\underline{u}}

\newcommand{\up}{\underline{p}}

\newcommand{\uf}{\underline{f}}
\newcommand{\ulam}{\underline{\lambda}}

\newcommand{\lap}{\Delta}
\newcommand{\laps}{\lap_{\M}}

\newcommand{\diag}{\text{diag}}
\DeclareMathOperator*{\argmin}{argmin}

\usepackage{alltt}
\usepackage{color}
\definecolor{string}{rgb}{0.7,0.0,0.0}
\definecolor{comment}{rgb}{0.13,0.54,0.13}
\definecolor{keyword}{rgb}{0.0,0.0,1.0}

\newcommand{\PAK}[1]{\textcolor{black}{#1}}

\definecolor{purple}{rgb}{0.5412,0.1686,0.8863}
\newcommand{\refa}[1]{{\color{black}#1}}
\newcommand{\refb}[1]{{\color{black}#1}}
\newcommand{\refc}[1]{{\color{black}#1}}

\begin{document}

\begin{frontmatter}

\author[1]{Andrew M.\ Jones}
\ead{andrewjones237@u.boisestate.edu}

\author[2]{Peter A.\ Bosler}
\ead{pabosle@sandia.gov}

\author[2]{Paul A.\ Kuberry}
\ead{pakuber@sandia.gov}

\author[1]{Grady B. Wright\corref{cor1}}
\ead{gradywright@boisestate.edu}

\cortext[cor1]{Corresponding author}

\affiliation[1]{organization={Boise State University}, addressline={1910 University Drive},
postcode={83725}, city={Boise}, state={Idaho}, country={USA}}

\affiliation[2]{organization={Sandia National Laboratories}, addressline={P.O. Box 5800},
postcode={87185}, city={Albuquerque}, state={New Mexico}, country={USA}}


%
\title{Generalized moving least squares vs.\ radial basis function finite difference methods for approximating surface derivatives}




\begin{abstract}
Approximating differential operators defined on two-dimensional surfaces is an important problem that arises in many areas of science and engineering.  Over the past ten years, localized meshfree methods based on generalized moving least squares (GMLS) and radial basis function finite differences (RBF-FD) have been shown to be effective for this task as they can give high orders of accuracy at low computational cost, and they can be applied to surfaces defined only by point clouds.  However, there have yet to be any studies that perform a direct comparison of these methods for approximating surface differential operators (SDOs).  The first purpose of this work is to fill that gap.  For this comparison, we focus on an RBF-FD method based on polyharmonic spline kernels and polynomials (PHS+Poly) since they are most closely related to the GMLS method. Additionally, we use a relatively new technique for approximating SDOs with RBF-FD called the tangent plane method since it is simpler than previous techniques and natural to use with PHS+Poly RBF-FD.  The second purpose of this work is to relate the tangent plane formulation of SDOs to the local coordinate formulation used in GMLS and to show that they are equivalent when the tangent space to the surface is known exactly.  The final purpose is to use ideas from the GMLS SDO formulation to derive a new RBF-FD method for approximating the tangent space for a point cloud surface when it is unknown.  For the numerical comparisons of the methods, we examine their convergence rates for approximating the surface gradient, divergence, and Laplacian as the point clouds are refined for various parameter choices. We also compare their efficiency in terms of accuracy per computational cost, both when including and excluding setup costs. 
\end{abstract}



\end{frontmatter}

\section{Introduction}\label{sec:intro}
The problem of approximating differential operators defined on two dimensional surfaces embedded in $\mathbb{R}^3$ arises in many multiphysics models.  For example, simulating atmospheric flows with Eulerian or Lagrangian numerical methods requires approximating the surface gradient, divergence, and Laplacian on the two-sphere~\cite{vallis_2006,AMJ:numericalweather07,AMJ:Fornberg15,AMJ:Bosler14}. Similar surface differential operators (SDOs) on more geometrically complex surfaces appear in models of ice sheet dynamics~\cite{Gowan2021}, biochemical signaling on cell membranes~\cite{Liue2104191118}, morphogenesis~\cite{stoop2015curvature}, texture synthesis~\cite{Mikkelsen2020Bump}, and sea-air hydrodynamics \cite{banaerjee2004surfdiv}.  

Localized meshfree methods based on generalized moving least squares (GMLS) and radial basis function finite differences (RBF-FD) have become increasingly popular over the last ten years for approximating SDOs and solving surface partial differential equations (PDEs); see, for example,~\cite{mahadevan2022metrics,LiangZhao13,TraskKuberry20,SUCHDE20192789,GrossEtAl20} for GMLS and~\cite{LSW2016,FlyerLehtoBlaiseWrightStCyr2012,FlyerWrightFornberg,SHANKAR2014JSC,PIRET2016,SHANKAR2018722,petras2018rbf,ShawThesis,GUNDERMAN2020109256,Wendland2020,Alvarez2021,wright2022mgm} for RBF-FD. \refb{These methods can be applied to surfaces defined by point clouds, without having to form a triangulation of the surface like surface finite element methods~\cite{dziuk_elliott_2013} or a level-set representation of the surface like embedded finite element methods~\cite{BertalmioEtAl2001}.   
Additionally, for the special case of the sphere, RBF-FD has been shown to be highly competitive with element based methods in terms of accuracy per degree of freedom~\cite{FlyerLehtoBlaiseWrightStCyr2012,GUNDERMAN2020109256,FlyerWrightFornberg}.} While there is one study dedicated to comparing GMLS and RBF-FD for approximating functions and derivatives in $\mathbb{R}^2$ and $\mathbb{R}^3$~\cite{Bayona2019}, there are no studies that compare them for approximating SDOs.  The first purpose of the present work is to fill this gap.


The RBF-FD methods referenced above use different approaches \PAK{for} approximating SDOs, while the GMLS methods essentially use the same approach based on (weighted) polynomial least squares. To keep the comparison to GMLS manageable, we will limit our focus to an RBF-FD method based on polyharmonic spline (PHS) kernels augmented with polynomials (or PHS+Poly) since they are most closely related to GMLS~\cite{Bayona2019}.  Additionally, these RBF-FD methods are becoming more and more prevalent as they can give high orders of accuracy that are controlled by the augmented polynomial degree~\cite{bayona2017role} and they do not require choosing a shape parameter, which can be computationally intensive to do in an automated way.  

The techniques for formulating SDOs also vary significantly in the RBF-FD methods referenced above, while the formulations used in GMLS are similar, being based on local coordinates to the surface.  In this work, we limit our focus to the so-called tangent plane formulation with RBF-FD, as it provides a more straightforward technique for incorporating polynomials in RBF-FD methods than~\cite{Alvarez2021,LSW2016,FlyerLehtoBlaiseWrightStCyr2012,SHANKAR2014JSC,PIRET2016,SHANKAR2018722,petras2018rbf,Wendland2020} and is related to the local coordinate formulation used in GMLS (see below).  Additionally, the comparison in~\cite{ShawThesis} of several RBF-FD methods for approximating the surface Laplacian (Laplace-Beltrami operator) revealed the tangent plane approach to be the most computationally efficient in terms of accuracy per computational cost.   The tangent plane method was first introduced by Demanet~\cite{DEMANET2006} for approximating the surface Laplacian using polynomial based approximations. Suchde \& Kuhnert~\cite{SUCHDE20192789} generalized this method to other SDOs using polynomial weighted least squares. Shaw~\cite{ShawThesis} (see also~\cite{wright2022mgm}) was the first to use this method for approximating the surface Laplacian with RBF-FD and Gunderman et.\ al.~\cite{GUNDERMAN2020109256} independently developed the method for RBF-FD specialized to the surface gradient and divergence on the unit two-sphere.  The second purpose of the present work is to analytically  compare the local coordinate formulation of SDOs used in GMLS to the tangent plane formulation and to show that these formulations are in fact identical when the tangent space for the surface is known exactly for the given point cloud.

When the tangent space is unknown, which is generally the case for surface represented by point clouds, it \refa{must be} approximated.  There has been little attention given in the literature on RBF-FD methods for how to do these approximations; commonly it is assumed that they are computed by some separate techniques (e.g.,~\cite{SHANKAR2014JSC,LSW2016,Wendland2020}).  However, for GMLS, these approximations are incorporated directly in the methods (e.g.,~\cite{LiangZhao13,TraskKuberry20,GrossEtAl20}).  The third purpose of this work is use the ideas from GMLS to develop a new RBF-FD technique for  approximating the tangent space directly using PHS+Poly.  By combining this with the tangent plane method, we arrive at the first comprehensive PHS+Poly RBF-FD framework for approximating SDOs on point cloud surfaces.


The GMLS and RBF-FD methods both use weighted combinations of function values over a local stencil of points to approximate SDOs.  They also feature a parameter $\ell$ for controlling the degree of polynomial precision of the formulas.  For the numerical comparisons of the methods, we investigate how the size of the stencils and the polynomial degree effect the convergence rates of the methods for approximating SDOs under refinement.  We focus on approximations of the surface gradient, divergence, and Laplacian operators on two topologically distinct surfaces, the unit two-sphere and the torus, which are representative of a broad range of application domains. In the case of the sphere, we also study the convergence rates of the methods for different point sets, including \refc{icosahedral  points that are popular in applications}.  Finally, we investigate the efficiency of the methods in terms of their accuracy versus computational cost, both when including and excluding setup costs.


Our numerical results demonstrate that RBF-FD and GMLS give similar convergence rates for the same choice of polynomial degree $\ell$, but overall RBF-FD results in lower errors.  We also show that the often-reported super convergence of GMLS for the surface Laplacian only happens for highly structured, quasi-uniform point sets, and when the point sets are more general (but still possibly quasi-uniform), this convergence rate drops to the theoretical rate.  Additionally, we find that the errors for RBF-FD can be further reduced with increasing stencil sizes, but that this does not generally hold for GMLS, and the errors can actually deteriorate.  Finally, we find that when setup costs are included, GMLS has an advantage in terms of efficiency, but if these are neglected then RBF-FD is more efficient.

The remainder of the paper is organized as follows. In Section \ref{sec:background}, we provide some background and notation on stencil-based approximations and on surface differential operators.  We follow this with a detailed overview of the GMLS and RBF methods in Section \ref{sec:gmls} and \ref{sec:rbffd_lap}, respectively.  In particular, Section \ref{sec:tangent_plane} shows the equivalence of the local coordinate and tangent plane formulations of some SDOs, and Section \ref{sec:rbffd_tangent_space} introduces an RBF-FD method for approximating the tangent space.  A comparison of some theoretical properties of the two methods is given in Section \ref{sec:theoretical}, while extensive numerical comparisons are given in Section \ref{sec:results}. We end with some concluding remarks in Section \ref{sec:remarks}.

\section{Background and notation\label{sec:background}}
\subsection{Stencils}\label{sec:stencils}
The RBF-FD and GMLS methods both discretize SDOs by weighted combinations of function values over a local \textit{stencil} of points.  This makes them similar to traditional finite-difference methods, but the lack of a grid, a tuple indexing scheme, and inherent awareness of neighboring points requires that some different notation and concepts be introduced. In this section we review the stencil notation that will be used in the subsequent sections.

Let $\X=\{\vx_i\}_{i=1}^N$ be a global set of points (point cloud) contained in some domain $\Omega$. A \textit{stencil of $\X$} is a subset of $n \leq N$ nodes of $\X$ that are close (see discussion below for what this means) to some point $\vx_{\rm c}\in\Omega$, which is called the \textit{stencil center}. In this work, the stencil center is some point from $\X$, so that $\vx_{\rm c} = \vx_i$, for some $1\leq i \leq N$, and this point is always included in the stencil.  We denote the subset of points making up the stencil with stencil center $\vx_i$ as $\X^i$ and allow the number of points in the stencil to vary with $\vx_i$.  To keep track \PAK{of} which points in $\X^i$ belong to $\X$, we use \textit{index set} notation and let $\sigma^i$ denote the set of indices of the $1 < n_i \leq N$ points from $\X$ that belong to $\X^i$.  Using this notation, we write the elements of the stencil as $\X^i = \{\vx_j\}_{j\in \sigma^i}$.  We also use the convention that the indices are sorted by the distance the stencil points are from the stencil center $\vx_i$, so that the first element of $\sigma^i$ is $i$.

With the above notation, we can define a general stencil-based approximation method to a given (scalar) linear differential operator $\mathcal{L}$.  Let $u$ be a scalar-valued function defined on $\Omega$ that is smooth enough so that $\mathcal{L}u$ is defined for all $\vx\in\Omega$.  The approximation to $\mathcal{L}u$ at any $\vx_i\in \X$ is given as
\begin{equation}\label{eq:stencils}
\mathcal{L}u|_{\vx = \vx_i} \approx \sum_{j\in\sigma^i} c_{ij}u(\vx_j).
\end{equation}
The weights $c_{ij}$ are determined by the method of approximation, which in this study will be either GMLS or RBF-FD.  These weights can be assembled into a sparse $N\times N$ ``stiffness'' matrix, similar to mesh-based methods.  Vector linear differential operators (e.g., the gradient) can be similarly defined where \eqref{eq:stencils} is used for each component and $\mathcal{L}$ is the scalar operator for that component.

\begin{figure}
\centering
\begin{tabular}{cc}
\includegraphics[width=0.4\textwidth]{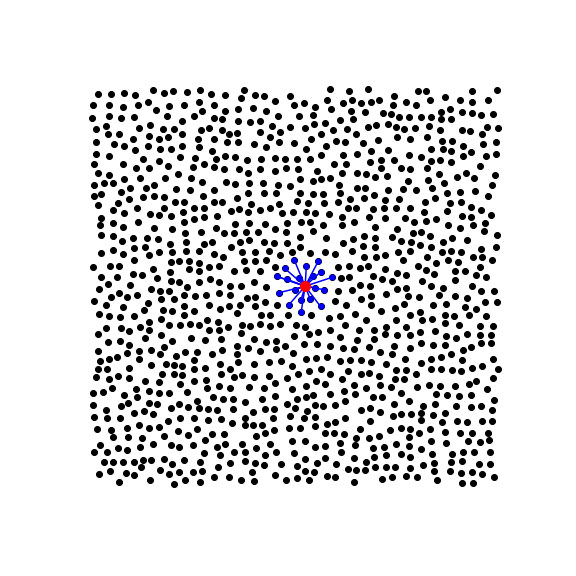} &
\includegraphics[width=0.4\textwidth]{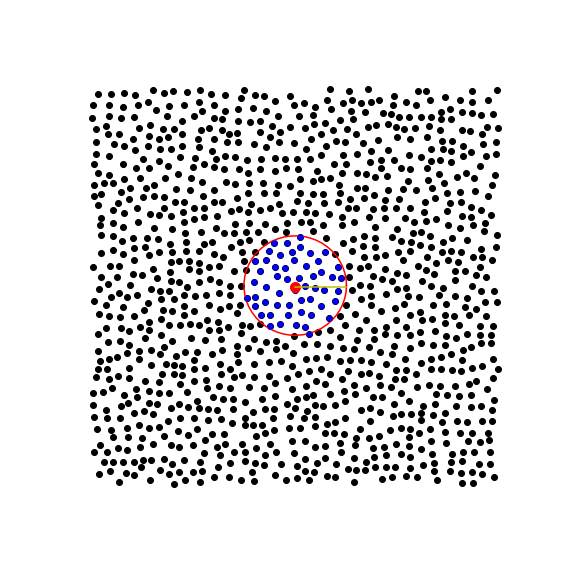} \\
(a) KNN & (b) Ball search
\end{tabular}
\caption{Comparison of the two search algorithms used in this paper for determining a stencil. The nodes $\X$ are marked with solid black disks and all the stencil points are marked with solid blue disks, except for the stencil center, which is marked in red.}\label{fig:search_algos}
\end{figure}

There are two main approaches used in the meshfree methods literature for determining the stencil points, one based on  $k$-nearest neighbors (KNN) and one based on ball searches.  These are illustrated in Figure \ref{fig:search_algos} for a scattered point set $\X$ in the plane.  The approach that uses KNN  is straightforward since it amounts to simply choosing the stencil $\X^i$ as the subset of $n_i$ points from $\X$ that are closest to $\vx_i$.  The approach that uses ball searches is a bit more involved, so we summarize it in Algorithm \ref{alg:epsilon_ball}.  Both methods attempt to select points such that the stencil satisfies polynomial unisolvency conditions (see the discussion in Section \ref{sec:metric_terms}).  In this work, we use the method in Algorithm \ref{alg:epsilon_ball} since
\begin{itemize}
\item it is better for producing stencils with symmetries when $\X$ is regular, which can be beneficial for improving the accuracy of the approximations;
\item it is more natural to use with the weighting kernel inherent to GMLS; and
\item it produces stencils that are not biased in one direction when the spacing of the points in $X$ are anisotropic.
\end{itemize}
To measure distance in the ball search, we use the standard Euclidean distance measured in $\mathbb{R}^3$ rather than distance on the surface since this is simple to compute for any surface.  We also use a $k$-d tree to efficiently implement \PAK{the} method.  Finally, the choice of parameters we use in Algorithm \ref{alg:epsilon_ball} are discussed in Section \ref{sec:gmls_stencils_weights}.  
\begin{algorithm}[t!]
\begin{algorithmic}[1]
\State \textbf{Input:} Point cloud $\X$; stencil center $\vx_{\rm c}$; number initial stencil points $n$; radius factor $\tau \geq 1$
\State \textbf{Output:} Indices $\sigma^{\rm c}$ in $\X$ for the stencil center $\vx_{\rm c}$
\State Find the $n$ nearest neighbors in $\X$ to $\vx_c$, using the Euclidean distance
\State Compute the max distance $h_{\max}$ between $\vx_{\rm c}$ and its $n$ nearest neighbors
\State Find the indices $\sigma^{\rm c}$ of the points in $\X$ contained in the ball of radius $\tau h_{\max}$ centered at $\vx_{\rm c}$
\caption{Procedure for determining the stencil points based on ball searches.\label{alg:epsilon_ball}}
\end{algorithmic}
\end{algorithm}


\subsection{Surface differential operators in local coordinates}\label{sec:local_coords}
Here we review some differential geometry concepts that will be used in the subsequent sections.  \refc{We refer the reader to the books~\cite{Walker15,ONeill2006,koenderink1990solid} for a thorough discussion of these concepts and the derivations of what follows.} 

\refc{We assume that $\M\subset\mathbb{R}^3$ is a regular surface and let $T_{\vx}\M$ denote the set of all vectors in $\mathbb{R}^3$ that are tangent to $\M$ at $\vx\in\M$ (i.e., the tangent space to $\M$ at $\vx$). This assumption implies that for each point $\vx\in\M$ there exists a local parameterization in $T_{\vx}\M$ of a neighborhood (or patch) of $\M$ containing $\vx$ of the form
\begin{align}
\vf(\xh,\yh) = (\xh,\yh,f(\xh,\yh)),
\label{eq:monge}
\end{align}
where $\xh$, $\yh$ are local coordinates for $T_{\vx}\M$, and $f$ is a smooth function for the ``height'' of the surface patch over $T_{\vx}\M$~\cite{koenderink1990solid}.  This local parametric representation of a surface is called a Monge patch or Monge form~\cite{ONeill2006} and is illustrated for a bumpy sphere surface in Figure \ref{fig:projection}.  As we see below, it is particularly well suited for computing SDOs.}

Using the parameterization \eqref{eq:monge}, the local metric tensor $G$ about $\vx$ for the surface is given as 
\begin{align}
G =
\begin{bmatrix}
\partial_{\xh} \vf \cdot \partial_{\xh} \vf & \partial_{\xh}\vf \cdot \partial_{\yh} \vf \\
\partial_{\yh} \vf \cdot \partial_{\xh} \vf & \partial_{\yh} \vf \cdot \partial_{\yh} \vf 
\end{bmatrix}=
\begin{bmatrix}
1 + (\partial_{\xh} f)^2 & (\partial_{\xh} f) (\partial_{\yh}f) \\
(\partial_{\xh} f) (\partial_{\yh}f) & 1 + (\partial_{\yh} f)^2
\end{bmatrix}.
\label{eq:metric_tensor}
\end{align}
Letting $g^{ij}$ denote the $(i,j)$ entry of $G^{-1}$, the surface gradient operator locally about $\vx$ is given as 
\begin{align}
\widehat{\nabla}_{\M} = (\partial_{\xh} \vf) \left(g^{11} \partial_{\xh} + g^{12} \partial_{\yh}\right) + 
(\partial_{\yh} \vf) \left(g^{21} \partial_{\xh} + g^{22} \partial_{\yh}\right).
\label{eq:intrinsic_grad_rot}
\end{align}
However, this is the surface gradient with respect to the horizontal $\xh\yh$-plane (see Figure \ref{fig:projection} (b)), and subsequently needs to be rotated so that it is with respect to $T_{\vx}\M$ in its original configuration.  If $\vxi^1$ and $\vxi^2$ are orthonormal vectors that span $T_{\vx}\M$ and $\veta$ is the unit outward normal to $\M$ at $\vx$, then the surface gradient in the correct orientation is given as
\begin{align}
{\nabla}_{\M}  = \underbrace{\begin{bmatrix} \vxi^1 & \vxi^2 & \veta\end{bmatrix}}_{\ds R} \widehat{\nabla}_{\M}.
\label{eq:intrinsic_grad}
\end{align}
Using this result, the surface divergence of a smooth vector $\vu\in T_{\vx}\M$ can be written as
\begin{align}
\nabla_{\M} \cdot \vu = \left(g^{11} \partial_{\xh} + g^{12} \partial_{\yh}\right) (\partial_{\xh} \vf)^T R^{T}\vu  + 
\left(g^{21} \partial_{\xh} + g^{22} \partial_{\yh}\right)(\partial_{\yh} \vf)^T R^{T}\vu
\label{eq:intrinsic_div}
\end{align}
\begin{figure}[htb]
\centering
\begin{tabular}{cc}
\includegraphics[width=0.4\textwidth]{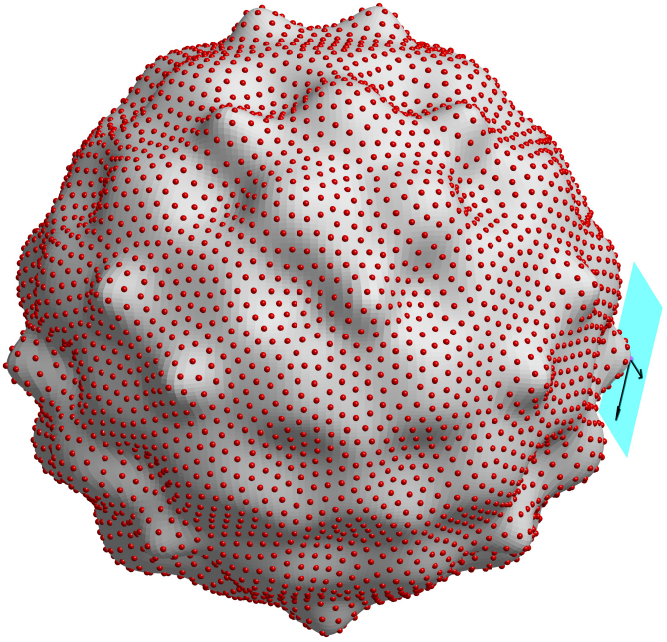} & \includegraphics[width=0.4\textwidth]{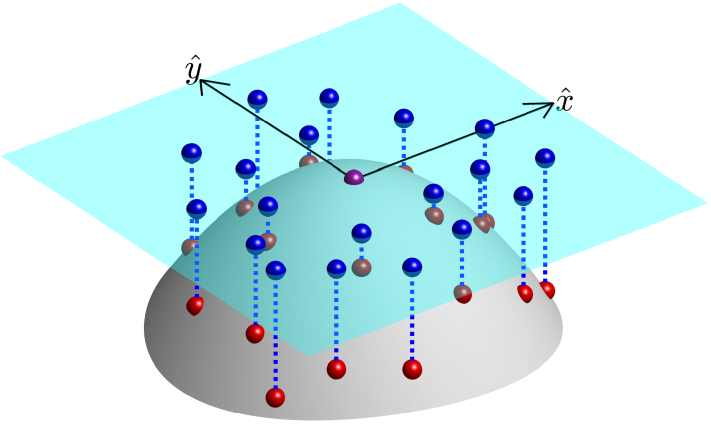} \\
(a) & (b) 
\end{tabular}
\caption{Illustration of a Monge patch parameterization for a local neighborhood of a regular surface $\M$ in 3D.  (a) Entire surface (in gray) together with the tangent plane (in cyan) for a point $\vx_{\rm c}$ where the Monge patch is constructed (i.e., $T_{\vx_{\rm c}}\M$); red spheres mark a global point cloud $\X$ on the surface.  (b) Close-up view of the Monge patch parameterization, together with the points from a stencil $\X_{\rm c}$ (red spheres) formed from $\X$ and the projection of the stencil to the tangent plane (blue spheres); the stencil center $\vx_{\rm c}$ is at the origin of the axes for the $\hat{x}\hat{y}$-plane and is marked with a violet sphere. \label{fig:projection}}
\end{figure}

The surface Laplacian operator locally about $\vx$ is given as
\begin{equation}
\begin{aligned}
\laps = \frac{1}{\sqrt{|g|}} \biggl(&
\partial_{\xh}\left(\sqrt{|g|}g^{11}\partial_{\xh}\right) + 
\partial_{\xh}\left(\sqrt{|g|}g^{12}\partial_{\yh}\right) + \\ &
\partial_{\yh}\left(\sqrt{|g|}g^{21}\partial_{\xh}\right) +
\partial_{\yh}\left(\sqrt{|g|}g^{22}\partial_{\yh}
\right)\biggr),
\end{aligned}
\label{eq:intrinsic_lap}
\end{equation}
where $|g|=\det(G)$. This operator is invariant to rotations of the surface in $\mathbb{R}^3$, so no subsequent modifications of \eqref{eq:intrinsic_lap} are necessary.

\section{GMLS using local coordinates}\label{sec:gmls}
The \PAK{formulation} of GMLS \PAK{on a manifold} was introduced by Liang \& Zhao~\cite{LiangZhao13} and further refined by Trask, Kuberry, and collaborators~\cite{TraskKuberry20,GrossEtAl20}.  It uses local coordinates to approximate SDOs as defined in \eqref{eq:intrinsic_grad}--\eqref{eq:intrinsic_lap} and requires a method to also approximate the metric terms. Both approximations are computed for each $\vx_i \in \X\subset \M$ using GMLS over a local stencil of points $\X^i\subset\X$.  Below we give a brief overview of the method assuming that the tangent/normal vectors for the surface are known for each $\vx_i\in\X$.  We then discuss a method for approximating these that is used in the Compadre Toolkit~\cite{AMJ:Compadre}, which we use in the numerical experiments.

We present the GMLS method through the lens of derivatives of MLS approximants as we feel this makes the analog to RBF-FD clearer, it is also closer to the description from~\cite{LiangZhao13}.  Other derivations of GMLS are based on weighted least squares approximants of general linear functionals given at some set of points, e.g.~\cite{Wendland:2004,Mirzaei11,trask2017high}.  However, both techniques produce the same result in the end~\cite{Mirzaei11}.  For a more thorough discussion of MLS approximants, see for example~\cite[ch. 22]{fasshauer2007meshfree} and the references therein.

\subsection{Approximating the metric terms\label{sec:metric_terms}}
The metric terms are approximated from an MLS reconstruction of the Monge patch of $\M$ centered at each target point $\vx_i$ using a local stencil of $n_i$ points $\X^i\subset\X$.  This procedure is illustrated in Figure \ref{fig:projection} and can be described as follows. First, the stencil $\X^i$ is expressed in the form of \eqref{eq:monge} (i.e., $(\xh_j,\yh_j,f_j)$, $j\in\sigma^i$), where $(\xh_j,\yh_j)$ are the coordinates for the stencil points in $T_{\vx_i}\M$, and $f_j = f(\xh_j,\yh_j)$ are samples of the surface as viewed from the $\xh\yh$-plane. These can be computed explicitly as
\begin{align}
\begin{bmatrix}
\xh_j \\
\yh_j \\
f_j
\end{bmatrix}
= {\underbrace{\begin{bmatrix} \vxi_i^1 & \vxi_i^2 & \veta_i\end{bmatrix}^T}_{\ds R_i^T}} (\vx_j - \vx_i),
\label{eq:projection}
\end{align}
where $\vxi_i^{1}$ and $\vxi_i^{2}$ are orthonormal vectors that span $T_{\vx_i}\M$ and $\veta_i$ is the unit normal to $\M$ at $\vx_i$.  To simplify the notation that follows, we let $\vxh_j=(\xh_j,\yh_j)$ and $\Xh^i = \{\vxh_j\}_{j\in\sigma^i}$ denote the projection of the stencil $\X^i$ to $T_{\vx_i}\M$.  Note that for convenience in what comes later we have shifted the coordinates so that the center of the projected stencil is $\vxh_i = (0,0)$.  

In the second step, the approximate Monge patch at $\vx_i$ is constructed from a MLS approximant of the data $(\vxh_j,f_j)$, $j\in\sigma^i$, which can be written as
\begin{align}
q(\vxh) = \sum_{k=1}^{L} b_k(\vxh) p_k(\vxh),
\label{eq:mls_interp}
\end{align}
where $\{p_1,\ldots,p_L\}$ is a basis for $\mathbb{P}^2_{\ell}$ (the space of bivariate polynomials of degree $\ell$) and $L = \text{dim}(\mathbb{P}^2_{\ell})=(\ell+1)(\ell+2)/2$ is the dimension of this space. The coefficients $b_k(\vxh)$ of the approximant are determined from the data according to the weighted least squares problem
\begin{align}
\ub^{*}(\vxh) = \argmin_{\ub\in\mathbb{R}^L} \sum_{j\in\sigma^i} w_{\rho}(\vxh_j,\vxh) (q(\vxh_j) - f_j)^2 =   \argmin_{\ub\in\mathbb{R}^L} \|W_{\rho}(\vxh)^{1/2}(P\ub - \uf)\|_2^2,
\label{eq:wls_fd_approx_system}
\end{align}
where $w_{\rho}:\mathbb{R}^2\times\mathbb{R}^2\rightarrow\mathbb{R}^{\geq 0}$ is a weight kernel that depends on a support parameter $\rho$,  $W_{\rho}(\vxh)$ is the $n_i\times n_i$ diagonal matrix $W_{\rho}(\vxh) = \diag(w_{\rho}(\vxh_j,\vxh))$, and $P$ is the $n_i\times L$ Vandermonde-type matrix
\begin{align}
P = \begin{bmatrix} p_1(\vxh_j) & p_2(\vxh_j) & \cdots & p_L(\vxh_j) \end{bmatrix},\; j\in\sigma^i
\label{eq:vandermonde}
\end{align}
Here we use underlines to denote vectors (i.e., $\ub$ and $\uf$ denote vectors containing coefficients and data from \eqref{eq:wls_fd_approx_system}, respectively).  Note that the coefficients $b_k$ depend on $\vxh$ because the kernel $w_{\rho}$ depends on $\vxh$ (this gives origin to the term ``moving'' in MLS). We discuss the selection of the stencils and weight\PAK{ing} kernel below, but for now it is assumed that $n_i > L$ and $\X^i$ is unisolvent on the space $\mathbb{P}^2_{\ell}$ (i.e., $P$ is full rank), so that \eqref{eq:wls_fd_approx_system} has a unique solution.  

The MLS approximant $q$ is used in place of $f$ in the Monge patch \eqref{eq:monge} and it is used to approximate the metric terms in \eqref{eq:intrinsic_grad}--\eqref{eq:intrinsic_lap}.  To compute these terms, various derivatives need to be approximated at the projected stencil center $\vxh_i$.  Considering, for example, $\partial_{\xh} q$, the approximation is computed as follows:
\begin{align}
\partial_{\xh} q\bigr|_{\vxh_i} \approx \sum_{k=1}^{L} b_k^{*}(\vxh_i) \partial_{\xh}(p_k(\vxh))\bigr|_{\vxh_i},
\label{eq:diff_mls_interp}
\end{align}
where $b_k^{*}(\vxh_i)$ come from \eqref{eq:wls_fd_approx_system} with $\vxh=\vxh_i$. Other derivatives of metric terms in \eqref{eq:intrinsic_grad}--\eqref{eq:intrinsic_lap} are approximated in a similar way to \eqref{eq:diff_mls_interp}.  We note if the standard monomial basis is used for $\{p_1,\ldots,p_L\}$, then by centering the projected stencil in \eqref{eq:projection} about the origin, only one of the derivatives of $p_k$ in \eqref{eq:diff_mls_interp} is non-zero when evaluated at $\vxh_i$. 

Note that \eqref{eq:diff_mls_interp} is only an approximation of $\partial_{\xh}q$ because it does not include the contribution of $\partial_{\xh}(b_k^{*}(\vxh))\big|_{\vxh_i}$.  This approximation is referred to as a ``diffuse derivative'' in the literature and is equivalent to the GMLS formulation of approximating derivatives~\cite{Mirzaei11}.  The term ``GMLS derivatives'' is preferred over ``diffuse derivatives'' to describe \eqref{eq:diff_mls_interp}, since the approximation is not diffuse or uncertain and has the same order of accuracy as the approximations that include the derivatives of the weight kernels~\cite{mirzaei2016error}.

\subsection{Approximating SDOs}\label{sec:gmls_sdos}
The procedure for approximating any of the SDOs in \eqref{eq:intrinsic_grad}--\eqref{eq:intrinsic_lap} is similar to the one for approximating the metric terms, but for this task we are interested in computing stencil weights as in \eqref{eq:stencils} instead of the value of a derivative at a point.  Since these SDOs involve computing various partial derivatives with respect to $\xh$ and $\yh$, we can use \eqref{eq:diff_mls_interp} as a starting point for generating these stencil weights.  If $\{u_j\}_{j\in\sigma^i}$ are samples of a function $u$ over the projected stencil $\Xh^i$, then we can again approximate $\partial_{\vxh}u\bigr|_{\vxh=\vxh_i}$ using \eqref{eq:diff_mls_interp}, with $b_k^*(\vxh_i)$ defined in terms of the samples of $u$.  To write this in stencil form we note that \eqref{eq:diff_mls_interp} can be written using vector inner products as
\begin{align}
\partial_{\xh} u\bigr|_{\vxh_i} \approx \partial_{\xh} q\bigr|_{\vxh_i} \approx {\underbrace{\begin{bmatrix} \partial_{\xh} p_1\bigr|_{\vxh_i} & \cdots & \partial_{\xh}p_L\bigr|_{\vxh_i}\end{bmatrix}}_{\ds (\partial_{\xh}\up(\vxh_i))^T}} \ub^{*}(\vxh_i) =
{\underbrace{\begin{bmatrix} c_{1}^{i} & \cdots & c_{n^{i}}\end{bmatrix}}_{\ds (\uc^i_{\xh})^T}} \uu,
\label{eq:wls_fd_approx}
\end{align}
where we have substituted the solution of $\ub^{*}(\vxh_i)$ in \eqref{eq:wls_fd_approx_system} to obtain the term in the last equality.  \refc{Using the normal equation solution for $\ub^{*}(\vxh_i)$, the  stencil weights $\uc^i_{\xh}$ can be expressed as
\begin{align}
\uc^i_{\xh} =  W_{\rho}(\vxh_i)P(P^T W_{\rho}(\vxh_i)P)^{-1} (\partial_{\xh}\up(\vxh_i)).
\label{eq:wls_fd_lap}
\end{align}
This is typically computed using a QR factorization of $W_{\rho}(\vxh_i)^{1/2}P$ to promote numerical stability.}

Stencil weights $\uc^i_{\yh}$, $\uc^i_{\xh\xh}$, $\uc^i_{\xh\yh}$, and $\uc^i_{\yh\yh}$ for the other derivative operators appearing in \eqref{eq:intrinsic_grad}--\eqref{eq:intrinsic_lap}
can be computed in a similar manner for each stencil $\Xh^i$, $i=1,\dots,N$.  These can then be combined together with the approximate metric terms to define the weights $\{c_{ij}\}$ in \eqref{eq:stencils} for any of the SDOs in  \eqref{eq:intrinsic_grad}--\eqref{eq:intrinsic_lap}.

\subsection{Choosing the stencils and weight kernel\label{sec:gmls_stencils_weights}}
As discussed in Section \ref{sec:stencils}, we use Algorithm \ref{alg:epsilon_ball} to choose the stencil weights.  For the initial stencil size, we use $L=\text{dim}(\mathbb{P}^2_{\ell})$.  The radius factor $\tau$ controls the size of the stencil, with larger $\tau$ resulting in larger stencils, and we experiment with this parameter in the numerical results section.

There are many choices for the weight kernel $w_{\rho}$ in \eqref{eq:wls_fd_approx_system}.  Typically, a single radial kernel is used to define $w_{\rho}$ as $w_{\rho}(\vx,\vy) = w(\|\vx-\vy\|/\rho)$, where $\vx,\vy\in\mathbb{R}^d$ and $\| \cdot \|$ is the standard Euclidean norm for $\mathbb{R}^d$.  In this work, we use the same family of compactly supported radial kernels as~\cite{TraskKuberry20,GrossEtAl20} and implemented in~\cite{AMJ:Compadre}:
\begin{align}
w_{\rho}(\vx,\vy) = \left(1 - \frac{\|\vx -\vy\|}{\rho}\right)_{+}^{2m},
\label{eq:weight_kernel}
\end{align}
where $m$ is a positive integer and $(\cdot)_{+}$ is the positive floor function.  These $C^{0}$ kernels have support over the ball of radius $\rho$ centered at $\vy$.  While smoother kernels can be used such as Gaussians, splines, or Wendland kernels~\cite{fasshauer2007meshfree}, we have not observed any significant improvement in the accuracy of GMLS derivative approximations with smoother kernels.  In general, proofs on how the choice of kernels effects the accuracy of GMLS approximations have yet to be found. 




Finally, we note that the support parameter $\rho$ is chosen on a per stencil basis and is set equal to $\tau h_{max}$ from Algorithm \ref{alg:epsilon_ball}.  \refb{Picking an optimal value for $\tau$ to minimize the approximation error is a difficult problem.  In general, the optimal value depends locally on the point set and the function (or its derivative) being approximated~\cite{lipman2006error}.  While there are some algorithms that attempt to approximate this value to minimize the local pointwise error (e.g.,~\cite{lipman2006error,wang2008optimal}), they are computationally expensive.  Typically, one chooses a single $\tau > 1$ such that the minimization problem \eqref{eq:wls_fd_approx_system} is well-posed (i.e., $P$ is full rank).  This can be easily monitored for each stencil to adjust $\tau$ appropriately.}


\subsection{Approximating the tangent space\label{sec:gmls_tangent_space}}
\begin{figure}[htb]
\centering
\begin{tabular}{cc}
\includegraphics[width=0.4\textwidth]{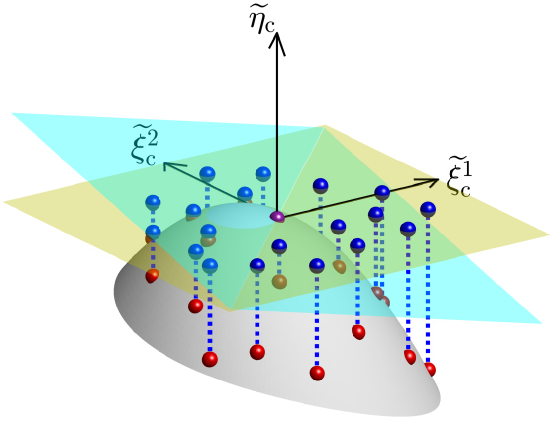} & \includegraphics[width=0.4\textwidth]{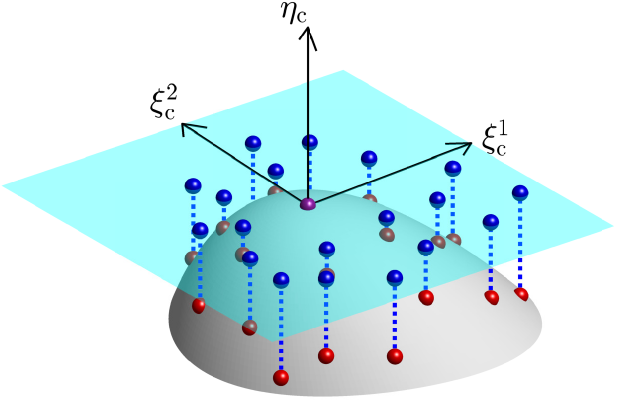} \\
(a) & (b) 
\end{tabular}
\caption{Illustration of the tangent plane correction method.  (a) Monge patch parameterization for a local neighborhood of a regular surface $\M$ (in gray) in 3D using a coarse approximation to the tangent plane (in yellow) at the center of the stencil $\vx_{\rm c}$ and the refined approximation to the tangent plane (in cyan).  (b) Same as (a), but for the Monge patch with respect to the refined tangent plane.  The red spheres denote the points from the stencil and the blue spheres mark the projection of the stencil to the (a) coarse and (b) refined tangent planes.  The coarse and refined approximations to the tangent and normal vectors are given as $\vxi^{1}_{\rm c}$, $\vxi^{2}_{\rm c}$, and $\veta_{\rm c}$, respectively, with tildes on these variables denoting the coarse approximation. \label{fig:tp_correction}}
\end{figure}

When the tangent space $T_{\vx_i}\M$ is unknown, a coarse approximation to it can be computed for each stencil $\X^i$ using principal component analysis~\cite{LiangZhao13}.  In this method, one computes the eigenvectors of the covariance matrix $\overline{X}_i \overline{X}_i^T$, where $\overline{X}_i$ is the $3$-by-$n_i$ matrix formed from the stencil points $\X^i$ centered about their mean.  The two dominant eigenvectors of this matrix are taken as a coarse approximation to $T_{\vx_{i}}\M$ and the third is taken as a coarse approximation to the normal to $\M$ at $\vx_i$; we denote these by $\tilde{\vxi}_{i}^1$, $\tilde{\vxi}_i^2$, and  $\tilde{\veta}_i$, respectively.  Next, an approximate Monge patch parameterization is formed with respect to this approximate tangent space using MLS following the same procedure outlined at the beginning of Section \ref{sec:metric_terms}.  This procedure is illustrated in Figure \ref{fig:tp_correction} (a), where the coarse approximate tangent plane is given in yellow.  A refined approximation to the true tangent plane and normal at the stencil center $\vx_{\rm c}$ can be obtained by computing the tangent plane and normal to the MLS approximant of the Monge patch at $\vx_{\rm c}$; this plane is given in cyan in Figure \ref{fig:tp_correction} (a).  Once this plane is computed, a new Monge patch parameterization with respect to this refined tangent plane approximation is formed, as illustrated in Figure \ref{fig:tp_correction} (b).  This procedure is repeated for each stencil $\X^i$ and the refined tangent space computed for each stencil is used in the procedure described in Section \ref{sec:metric_terms} for approximating the metric terms. 

\section{RBF-FD using the tangent plane}\label{sec:rbffd_lap}
As discussed in the introduction, there are several RBF-FD methods that have been developed over the past ten years for approximating SDOs.  We use the one based on the tangent plane method for formulating SDOs and PHS+Poly interpolants for approximating the derivatives that appear in this formulation.  The subsections below provide a detailed overview of these respective techniques.

\subsection{Tangent plane method\label{sec:tangent_plane}}
The tangent plane method similarly uses local coordinates for the surface in the tangent plane formed at each $\vx_i\in \X$, but unlike the method from Section \ref{sec:metric_terms}, it does not use approximations to the metric terms.  It instead approximates the SDOs at each $\vx_i$ using the standard definitions for the derivatives in the tangent plane.  So, using local coordinates \eqref{eq:monge} about $\vx_i$, \refa{the} surface gradient for the tangent plane method is taken as
\begin{align}
{\nabla}_{\M}  = R_i \begin{bmatrix} \partial_{\xh} \\ \partial_{\yh} \\  0 \end{bmatrix},
\label{eq:tp_grad}
\end{align}
and the surface divergence of a smooth vector $\vu\in T_{\vx_i}\M$ is taken as
\begin{align}
\nabla_{\M} \cdot \vu = 
\begin{bmatrix} \partial_{\xh} & \partial_{\yh}  & 0 \end{bmatrix}R_i^{T}\vu,
\label{eq:tp_div}
\end{align}
where $R_i$ is the rotation matrix given in \eqref{eq:projection}.  Similarly, the surface Laplacian in the tangent plane method is 
\begin{equation}
\laps = \partial_{\xh\xh} + \partial_{\yh\yh}.
\label{eq:tp_lap}
\end{equation}
We next show that if $T_{\vx_i}\M$ is known exactly for each $\vx_i\in \X$ and the point at which the SDOs are evaluated is $\vx_i$, then the SDOs \eqref{eq:tp_grad}--\eqref{eq:tp_lap} are equivalent to the corresponding ones involving metric terms \eqref{eq:intrinsic_grad}--\eqref{eq:intrinsic_lap}.  This was shown indirectly in~\cite{DEMANET2006} for the surface Laplacian using the distributional definition of the surface Laplacian.  Here we show the result follows explicitly for each surface differential operator \eqref{eq:intrinsic_grad}--\eqref{eq:intrinsic_lap} from the local coordinate formulation in Section \ref{sec:local_coords}.


The first step is to note that the vectors $\partial_{\xh} \vf\bigr|_{\vxh_i}$ and $\partial_{\yh} \vf\bigr|_{\vxh_i}$ from the Monge parameterization \eqref{eq:monge} are tangential to the $\xh\yh$-plane and must therefore be orthogonal to the vector $\begin{bmatrix} 0 & 0 & 1 \end{bmatrix}$.  This implies $\partial_{\xh} f= \partial_{\yh}f = 0$ at $\vxh_i$, which means the metric tensor \eqref{eq:metric_tensor} reduces to the identity matrix when evaluated at $\vxh_i$.  Using this result in \eqref{eq:intrinsic_grad} for $\widehat{\nabla}_{\M}$ means that the surface gradient formula \eqref{eq:intrinsic_grad} is exactly \eqref{eq:tp_grad} when evaluated at $\vxh_i$.  The equivalence of the surface divergence formulas \eqref{eq:intrinsic_div} and \eqref{eq:tp_div} also follow immediately from this result.

The steps for showing the equivalence of the surface Laplacian operator are more involved. To simplify the notation in showing this result, we denote partial derivatives of $f$ with subscripts.  For the first step of this process, we substitute the explicit metric terms, $|g| = (1+f_{\xh}^2)(1+f_{\yh}^2) - (f_{\xh} f_{\yh})^2$, $g^{11} = (1+f_{\yh})/|g|$, $g^{12}=g^{21}=-(f_{\xh})(f_{\yh})/|g|$, and $g^{22} = (1+f_{\xh})/|g|$, into \eqref{eq:intrinsic_lap} and expand the derivatives.  Next, we simplify to obtain the following formula:
\begin{align*}
\begin{aligned}
\laps = \frac{1}{\left(f_{\xh}^2+f_{\yh}^2+1\right)^{2}} \biggl( &
\left(f_{\yh} f_{\xh\yh} \left(1+2 f_{\xh}^2+f_{\yh}^2\right)-(f_{\xh} f_{\xh\xh} + f_{\yh}f_{\xh\yh})(1+f_{\yh}^2) - f_{\xh}f_{\yh\yh}(1 + f_{\xh}^2)\right)\partial_{\xh} + \\ 
& \left(f_{\xh} f_{\xh\yh} \left(1+2 f_{\yh}^2+f_{\xh}^2\right)-(f_{\yh} f_{\yh\yh}  + f_{\xh}f_{\xh\yh})(1 + f_{\xh}^2) - f_{\yh}f_{\xh\xh}(1+f_{\yh}^2) \right)\partial_{\yh}\biggr) + \\
& g^{11} \partial_{\xh\xh} + 2g^{12} \partial_{\xh\yh} + g^{22}\partial_{\yh\yh}
\end{aligned}
\end{align*}
Using $f_{\xh} = f_{\yh} = g^{12} = 0$ and $g^{11}=g^{22}=1$ at $\vxh_i$, this formula reduces to \eqref{eq:tp_lap}.

\subsection{Approximating the SDOs}
Since the tangent plane method does not require computing approximations to any metric terms, we only need to describe the RBF-FD method for approximating the derivatives that appear in \eqref{eq:tp_grad}--\eqref{eq:tp_lap}.  We derive this method from derivatives of interpolants over the projected stencils for each point $\vx_i\in X$ using the same notation as Section \ref{sec:gmls} and we assume that the tangent space is known.  A method for approximating the tangent space also using RBF-FD is discussed in Section \ref{sec:rbffd_tangent_space}.



Let $\{u_j\}_{j\in\sigma^i}$ be samples of some function $u$ over the projected stencil $\Xh^i = \{\vxh_j\}_{j\in\sigma^i}$.  The 
PHS+Poly interpolant to this data can be written
\begin{align}
s(\vxh) = \sum_{j=1}^{n_i} a_j \phi(\|\vxh - \vxh_{\sigma^i_j}\|) + \sum_{k=1}^{L} b_k p_k(\vxh),
\label{eq:phs_interp}
\end{align}
where $\phi(r) = r^{2\kappa+1}$ is the PHS kernel of order $2\kappa+1$, $\kappa\in\mathbb{Z}^{\geq 0}$, $\sigma^i_j$ is the $j$th index in $\sigma^i$, $\|\cdot\|$ denotes the Euclidean norm, and $\{p_1,\ldots,p_L\}$ are a basis for $\mathbb{P}^2_{\ell}$. The expansion coefficients are determined by the $n_i$ interpolation conditions and $L$ additional moment conditions:
\begin{align}
s(\vxh_j) = u_j,\; j\in\sigma^i \quad \text{and} \quad \sum_{j=1}^{n_i} a_j p_k(\vxh_{\sigma^i_j}) = 0,\; k=1,\ldots,L.
\label{eq:phs_interp_conditions}
\end{align}
These conditions can be written as the following $(n_i+L)\times(n_i+L)$ linear system
\begin{align}
\begin{bmatrix}
A & P \\
P^T & \vect{0}
\end{bmatrix}
\begin{bmatrix}
\ua \\
\ub
\end{bmatrix}
=
\begin{bmatrix}
\uu \\
\underline{0}
\end{bmatrix},
\label{eq:rbf_fd_interp_system}
\end{align}
where $A_{jk} = \|\vxh_{\sigma_j^i} - \vxh_{\sigma_k^i}\|^{2\kappa+1}$ ($j,k=1,\ldots,n_i$) and $P$ is the same Vandermonde-type matrix given in \eqref{eq:vandermonde}.  The PHS  parameter $\kappa$ controls the smoothness of the kernel and should be chosen such that $0\leq \kappa \leq \ell$.  With this restriction on $\kappa$, it can be shown that $A$ is positive definite on the subspace of vectors in $\mathbb{R}^n$ satisfying the $L$ moment conditions in \eqref{eq:phs_interp_conditions}~\cite{Wendland:2004}.  Hence, if the stencil points $\X^i$ are such that $\text{rank}(P)=L$ (i.e., they are unisolvent on the space $\mathbb{P}^2_{\ell}$), then the system \eqref{eq:rbf_fd_interp_system} is non-singular and the PHS+Poly interpolant is well-posed.  Note that this is the same restriction on $\X^i$ for the MLS problem \eqref{eq:wls_fd_approx_system} to have a unique solution.  

The stencil weights for approximating any of the derivatives appearing in the SDOs \eqref{eq:tp_grad}-\eqref{eq:tp_lap} can be obtained from differentiating the PHS+Poly interpolant \eqref{eq:phs_interp}.   Without loss of generality, consider approximating the operator $\partial_{\xh}$ over the stencil $\Xh^i$.  Using vector inner products as in \eqref{eq:wls_fd_approx},  the stencil weights for this operator are determined from the approximation
\begin{align*}
\partial_{\xh} u\bigr|_{\vxh_i} \approx \partial_{\xh} s\bigr|_{\vxh_i} 
&= \begin{bmatrix} \partial_{\xh}\underline{\phi}(\vxh_i) & \partial_{\xh}\up(\vxh_i) \end{bmatrix}^T\begin{bmatrix} \ua \\ \ub \end{bmatrix}.
\end{align*}
where $\partial_{\xh}\underline{\phi}(\vxh_i)$  and $\partial_{\xh}\up(\vxh_i)$ are vectors containing the entries $\partial_{\xh}\|\vxh - \vxh_{\sigma^i_j}\|^{2\kappa+1}\bigr|_{\vxh_i}$, $j=1,\ldots,n_i$, and $\partial_{\xh}p_k(\vxh)\bigr|_{\vxh_i}$, $k=1,\ldots,L$, respectively. Using \eqref{eq:rbf_fd_interp_system} in the preceding expression gives the stencil weights as the solution to the following linear system
\begin{align}
\begin{bmatrix}
A & P \\
P^T & \vect{0}
\end{bmatrix}
\begin{bmatrix}
\uc_{\xh}^i \\
\ulam
\end{bmatrix}
=
\begin{bmatrix}
\partial_{\xh}\underline{\phi}(\vxh_i) \\
\partial_{\xh}\up(\vxh_i)
\end{bmatrix},
\label{eq:rbf_fd_lap}
\end{align}
where the entries in $\ulam$ are not used as part of the weights.  \refc{Note that this description is equivalent to applying $\partial_{\xh}$ to the PHS+Poly cardinal basis functions defined over the stencil and evaluating them at $\vxh_i$~\cite{AMJ:Fornberg15}.}

Stencil weights $\uc^i_{\yh}$, $\uc^i_{\xh\xh}$, $\uc^i_{\xh\yh}$, and $\uc^i_{\yh\yh}$ for the other partial derivatives can be computed in an analogous way for each stencil $\Xh^i$, $i=1,\dots,N$.  These can then be combined together to define the weights $\{c_{ij}\}$ in \eqref{eq:stencils} for any of the SDOs in  \eqref{eq:tp_grad}--\eqref{eq:tp_lap}.

\subsection{Choosing the stencils and PHS order\label{sec:rbffd_stencils_weights}}
Similar to GMLS, we use Algorithm \ref{alg:epsilon_ball} to choose the stencils and also use the same initial stencil size of $n=L$ for this algorithm.  The parameter $\kappa$ used to determine the PHS order should be chosen with an upper bound of $\kappa \leq \ell$ (so that \eqref{eq:rbf_fd_lap} is well posed) and a lower bound such that the derivatives of the PHS kernels make sense for whatever operator the RBF-FD stencils are being used to approximate.  In this work we use $\kappa=\ell$ as \PAK{we} have found \PAK{that} this choice works well for approximating various SDOs across a wide range of surfaces.  Choosing $\kappa < \ell$ can be useful for improving the conditioning of the system \eqref{eq:rbf_fd_lap} and for reducing Runge Phenomenon-type edge effects in RBF-FD approximations near boundaries~\cite{bayona2019role}.


\subsection{Approximating the tangent space\label{sec:rbffd_tangent_space}}
If $T_{\vx_i}\M$ is unknown for any $\vx_i\in \X$, then we use a similar procedure to the one discussed for GMLS in Section \ref{sec:gmls_tangent_space} (and illustrated in Figure \ref{fig:tp_correction}) to approximate it.  The difference for RBF-FD is that instead of using an MLS reconstruction of the Monge patch parameterization formed from the coarse tangent plane approximation at each $\vx_i$, we use the PHS+Poly interpolant \eqref{eq:phs_interp} for the reconstruction.  The refined approximation to the tangent plane at each $\vx_i$ is then obtained from derivatives of the PHS+Poly interpolant of the Monge patch for stencil $\X^i$.  We note that this approach is new amongst the different tangent plane methods, as previous approaches assumed the tangent space was computed by some other, possibly unrelated techniques, and not directly from the stencils (e.g.,~\cite{SUCHDE20192789,ShawThesis,wright2022mgm}).  By combining this technique with the tangent plane method, we arrive at the first comprehensive PHS+Poly RBF-FD framework for approximating SDOs on point cloud surfaces.


\section{Theoretical comparison of GMLS and RBF-FD\label{sec:theoretical}}
In this section, we make comparisons of the GMLS and RBF-FD methods in terms of some of their theoretical properties, including the different approaches in formulating SDOs, the parameters of the approximations, and the computational cost.

One of the main differences between the GMLS and RBF-FD approaches is that the former uses the local coordinate method to formulate SDOs, while the latter uses the tangent plane method.  As shown in Section \ref{sec:tangent_plane} these methods are equivalent if the tangent space for $\M$ is known for each $\vx_i\in X$ and the SDOs are evaluated at the stencil center $\vx_i$. However, the GMLS method does not take advantage of this and instead includes metric terms in the formulation.  These metric terms are approximated with the same order of accuracy as the GMLS approximation of the derivatives (see below), so that these errors are asymptotically equivalent as the spacing of the points in the stencil goes to zero.  When the tangent space is unknown, both methods again approximate it to the same order of accuracy as their respective approximations of the derivatives.  

 

The GMLS and RBF-FD methods each feature the parameter $\ell$, which controls the degree of the polynomials used in the approximation.  For a given $\ell$, the formulas for either method are exact for all bivariate polynomials of degree $\ell$ in the tangent plane formed by the stencil center $\vx_i$.  Unsurprisingly, $\ell$ also effects the local accuracy of the formulas in the tangent plane with increasing $\ell$ giving higher orders of accuracy for smooth problems; see~\cite{mirzaei2016error,LiangZhao13} for a study of the accuracy of GMLS and~\cite{davydov2019optimal,bayona2019insight} for RBF-FD.  The order of accuracy of both methods depends on the highest order derivative appearing in the SDOs, and is generally $\ell$ if th\PAK{e derivative order} is 1 and $\ell-1$ if th\PAK{e derivative} order is two.  However, for certain quasi-uniform point clouds with symmetries, the order has been shown to be $\ell$ for GMLS applied to second order operators like the surface Laplacian~\cite{LiangZhao13}.


The computational cost of the methods can be split between the setup cost and the evaluation cost.  The setup cost depends on $\ell$ and $n_i$ (which depends on $\tau$).  For each stencil $\X^i$, the dominant setup cost of GMLS comes from solving the $n_i \times L$  system \eqref{eq:wls_fd_lap}, while the dominant cost for RBF-FD comes from solving the $(n_i + L)\times(n_i+L)$ system \eqref{eq:rbf_fd_lap}.  We use QR factorization to solve the GMLS system and LU factorization to solve the RBF-FD system, which gives the following (to leading order):
\begin{align}
\text{Setup cost GMLS} \sim 2 \sum_{i=1}^{N} n_i L^2\quad\text{and}\quad \text{Setup cost RBF-FD} \sim \frac23 \sum_{i=1}^{N} (n_i + L)^3.
\label{eq:setup_cost}
\end{align}
The stencil sizes depend on $\ell$ and $\tau$, and for quasi-uniform point clouds $\X$, $n_i$ is typically some multiple $\gamma$ of $L$.  In this case, the setup cost of RBF-FD is higher by approximately $(1+\gamma)^3/(3\gamma)$. We note that the setup procedures for both \refa{methods are} an embarrassingly parallel process, as each set of stencil weights can be computed independently of every other set.  The evaluation costs of both methods are the same and can be reduced to doing sparse matrix-vector products.  So, for a scalar SDO like the surface Laplacian
\begin{align}
\text{evaluation cost GMLS \& RBF-FD:}\sim 2\sum_{i=1}^N n_i.
\label{eq:runtime_cost}
\end{align}
If the $\ell$ and $\tau$ parameters remain fixed so that size of the stencils remain fixed as $N$ increases, then both the setup and evaluation cost are linear in $N$.

\section{Numerical comparison of GMLS and RBF-FD}\label{sec:results}
\begin{figure}[t]
\centering
\begin{tabular}{ccc}
\includegraphics[width=0.28\textwidth]{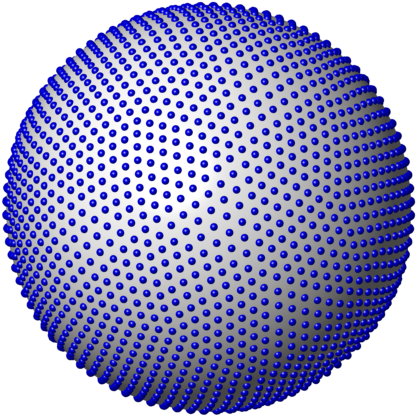} &
\includegraphics[width=0.28\textwidth]{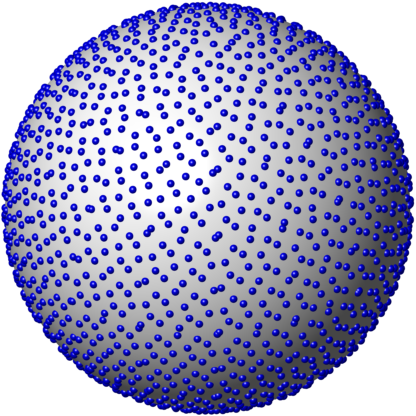} &
\includegraphics[width=0.33\textwidth]{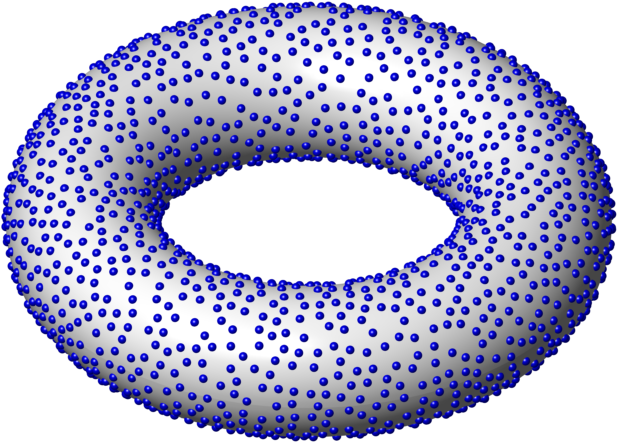} \\
\refa{(a) Icosahedral} & \refa{(b) Hammersley} & \refa{(c) Poisson disk}
\end{tabular}
\caption{\refa{Examples from the three node sets considered in the numerical experiments: (a) $N=2562$, (b) $N=2048$, (c) $N=2038$.}\label{fig:node_sets}}
\end{figure}

We perform a number of numerical experiments comparing GMLS and RBF-FD for approximating the gradient, divergence, and Laplacian on two topologically distinct surfaces: the unit two sphere $\Sph^2$ and the torus defined implicitly as
\begin{align}
\mathbb{T}^2 = \left\{(x,y,z)\in\mathbb{R}^3\, \bigr|\, (1-\sqrt{x^2 + y^2})^2 + z^2 - 1/9 = 0\right\}. \label{eq:torus}
\end{align}
\refa{For the experiments with the sphere, we consider two different node sets $\X$, icosahedral and Hammersley; see Figure \ref{fig:node_sets} (a) \& (b) for examples.  The first are highly structured, quasi-uniform points that are commonly used in numerical weather prediction~\cite{AMJ:Bosler14,AMJ:numericalweather07}. They have also been used in other studies on GMLS~\cite{TraskKuberry20} and RBF-FD~\cite{FlyerLehtoBlaiseWrightStCyr2012} methods on the sphere.  Hammersely are low discrepancy point sequences commonly used in Monte-Carlo integration on the sphere~\cite{Cui97equidistribution}.  They are highly unstructured with some points that nearly overlap.  For the experiments on the torus, we use Poisson disk points generated using the weighted sample elimination (WSE) algorithm~\cite{Yuksel2015}.  These points are also unstructured, but are quasi-uniform; see Figure \ref{fig:node_sets} (c) for an example.  They have also previously been used in studies on GMLS and RBF-FD methods~\cite{wright2022mgm}.}  \refc{Convergence results with other point sets can be found in the PhD thesis of the first author~\cite{JonesThesis}.}  

\refa{Error estimates for GMLS and RBF-FD typically require the nodes to be quasi-uniform in the sense that the average spacing between the points $h$ (or more generally the mesh-norm) decreases like $h\sim N^{-1/2}$~\cite{fasshauer2007meshfree,Wendland:2004}.  As mentioned above, the icosahedral and Poisson disk node sets have this property and are thus well-suited for numerically testing convergence rates of GMLS and RBF-FD methods with increasing $N$ (i.e., convergence as the density of the sampling of the surfaces increases).  Specifically, we experimentally examine the algebraic convergence rates $\beta$ versus $\sqrt{N}$, assuming the error behaves like $\mathcal{O}(N^{-\beta/2})$, and include results for polynomial degrees $\ell=2$, 4, and 6.  The Hammersley node sets are well suited to testing how stable the methods are to stencils with badly placed points.  Since these nodes have low discrepancy over the sphere, it also makes sense to test convergence in a similar manner to the other point sets.  The exact values of $N$ used in the experiments for each of the node sets are as follows.  For Icosahedral, $N=10242, 40962, 163842, 655362$, and for Hammersley and Poisson disk: $N=8153, 32615, 130463, 521855$.}




All RBF-FD results that follow were obtained from a Python implementation of the method that only utilizes the scientific computing libraries SciPy and NumPy. For the GMLS results, we use the software package Compadre~\cite{AMJ:Compadre}, which is implemented in C++ and uses the portable performance library Kokkos. 
\subsection{Convergence comparison: Sphere}\label{subsec:sphere}
We base all the convergence comparisons for the sphere on the following function consisting of a random linear combination of translates of $50$ Gaussians of different widths on the sphere: 
\begin{equation}
u(\vx) = \sum_{j=1}^{50} d_{j} \exp(-\gamma_j \|\vx- \vy_{j}\|^2),\; \vx,\vy_j\in\Sph^2,
\label{eq:target_sphere}
\end{equation}
where $\vy_{j}$ are the centers and are randomly placed on the sphere, and $d_j$ \& $\gamma_j$ are sampled from the normal distributions $\mathcal{N}(0,1)$ \& $\mathcal{N}(15,4)$, respectively.  This function has also been used in other studies on RBF-FD methods~\cite{LSW2016}.  We use samples of $u$ in the surface gradient tests and measure the error against the exact surface gradient, which can be computed using the Cartesian gradient $\nabla$ in $\mathbb{R}^3$ as $\nabla_\M u = \nabla u - \veta(\veta\cdot \nabla u)$, where $\veta$ is the unit outward normal to $\Sph^2$~\cite{FuselierWright2013} (which is just $\vx$).  Applying this to \eqref{eq:target_sphere} gives
\begin{align}
\nabla_{\M} u = 2\sum_{j=1}^{50} d_{j} \gamma_j (\vy_j - \vx(\vx\cdot\vy_j))\exp(-\gamma_j \|\vx- \vy_{j}\|^2).
\label{eq:grad_sphere}
\end{align}
We use samples of this field in the surface divergence tests.  Since $\nabla_{\M}\cdot\nabla_{\M} u = \laps u$, we compare the errors in this test against the exact surface Laplacian of $u$, which can be computed using the results of~\cite{AMJ:Fornberg15} as
\begin{align*}
\laps u = -\sum_{j=1}^{50} d_{j} \gamma_j (4 - \|\vx- \vy_{j}\|^2(2 + \gamma_j (4 - \|\vx- \vy_{j}\|^2)))\exp(-\gamma_j \|\vx- \vy_{j}\|^2).
\end{align*}
We also use this in the tests of the surface Laplacian using samples of $u$.   

For all these tests, we set radius factor $\tau$ in the stencil selection Algorithm \ref{alg:epsilon_ball} to 1.5, \refa{which gave good results for both RBF-FD and GMLS} (see the next section for some results on the effects of increasing $\tau$).  While the exact tangent space for the sphere is trivially determined, we approximate it in all the results using the methods discussed in the Section \ref{sec:gmls_tangent_space} for GMLS and Section \ref{sec:rbffd_tangent_space} for RBF-FD.  These approximations are done with the same parameters for approximating the different SDOs to keep the asymptotic orders of accuracy comparable.  Although not included here, we did experiments with the exact tangent space and obtained similar results to those presented here.

\begin{figure}
\centering
\begin{tabular}{ccc}
\rotatebox{90}{\hspace{0.1\textwidth}\small \quad \textbf{(a) Gradient}} & \includegraphics[width=0.4\textwidth]{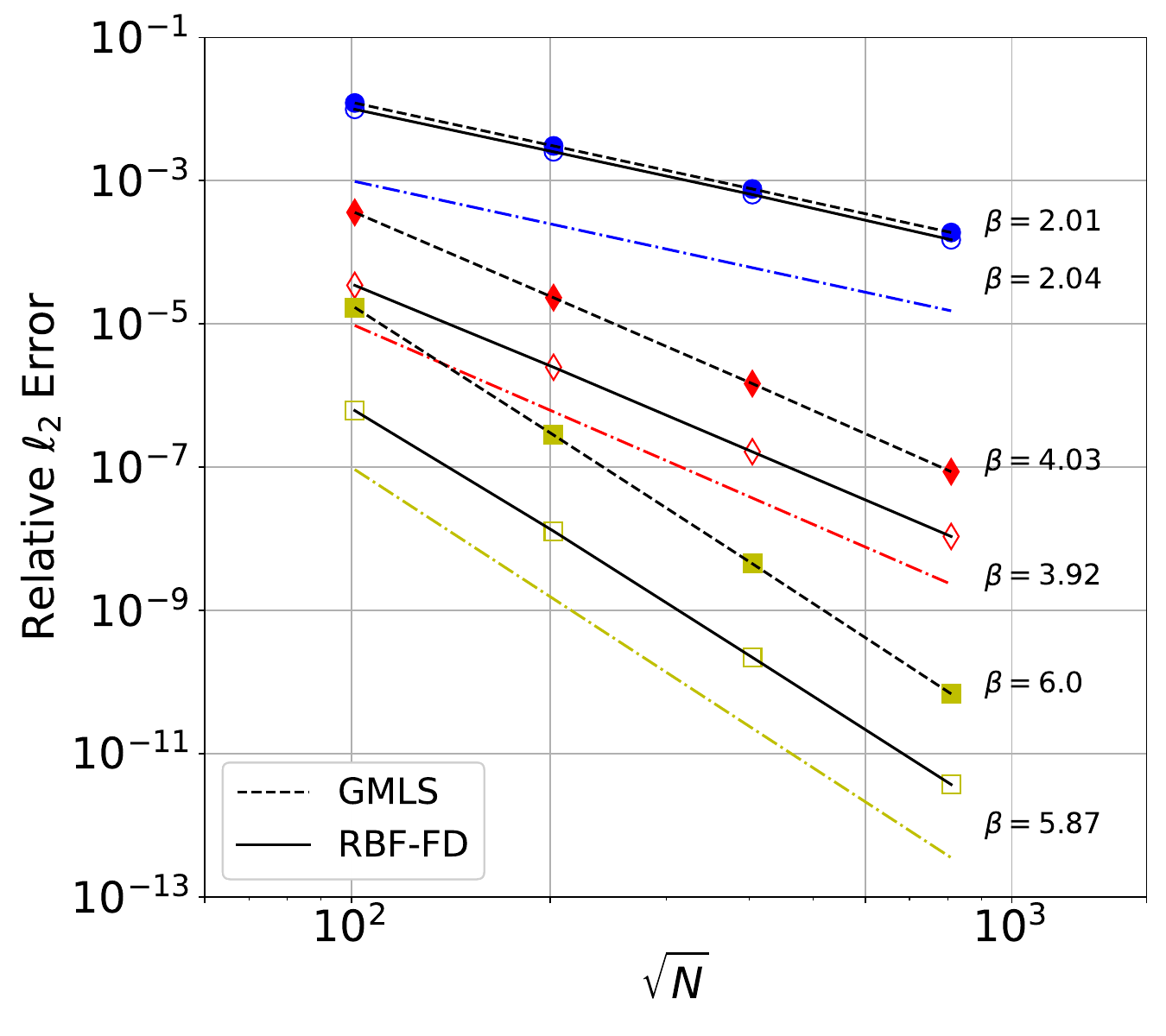} &
\includegraphics[width=0.4\textwidth]{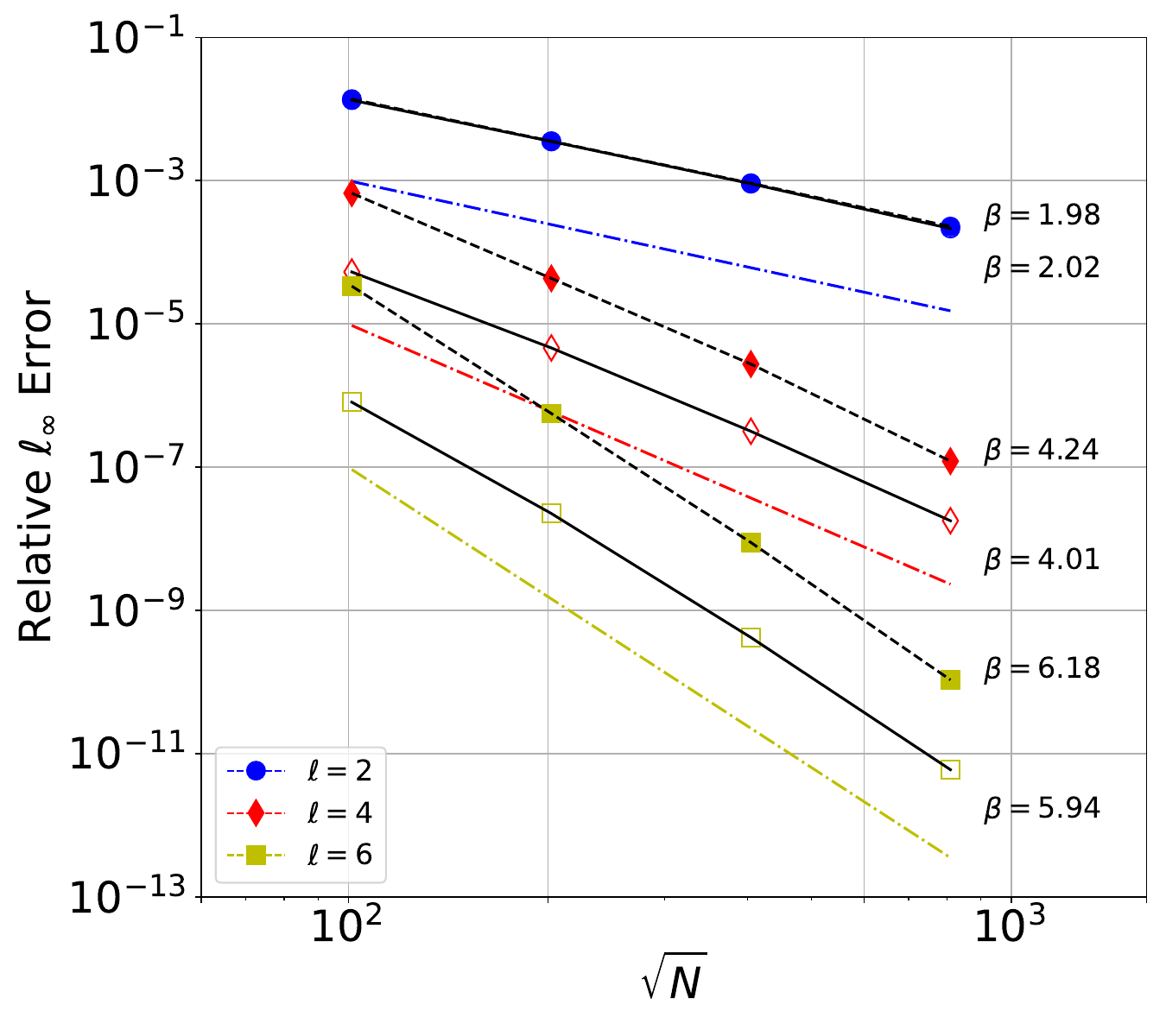} \\
\rotatebox{90}{\hspace{0.1\textwidth} \small \quad \textbf{(b) Divergence}} & \includegraphics[width=0.4\textwidth]{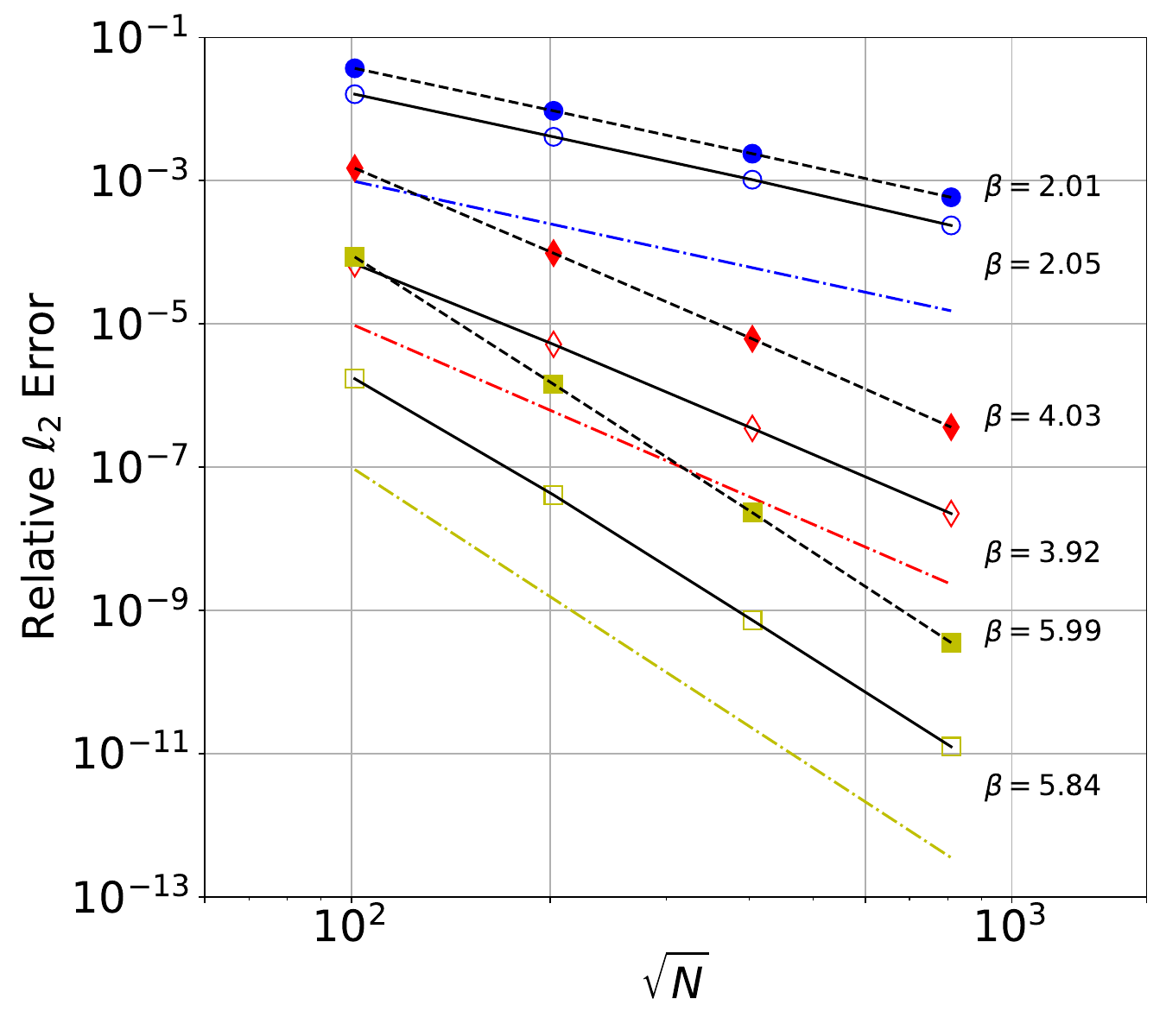} & 
\includegraphics[width=0.4\textwidth]{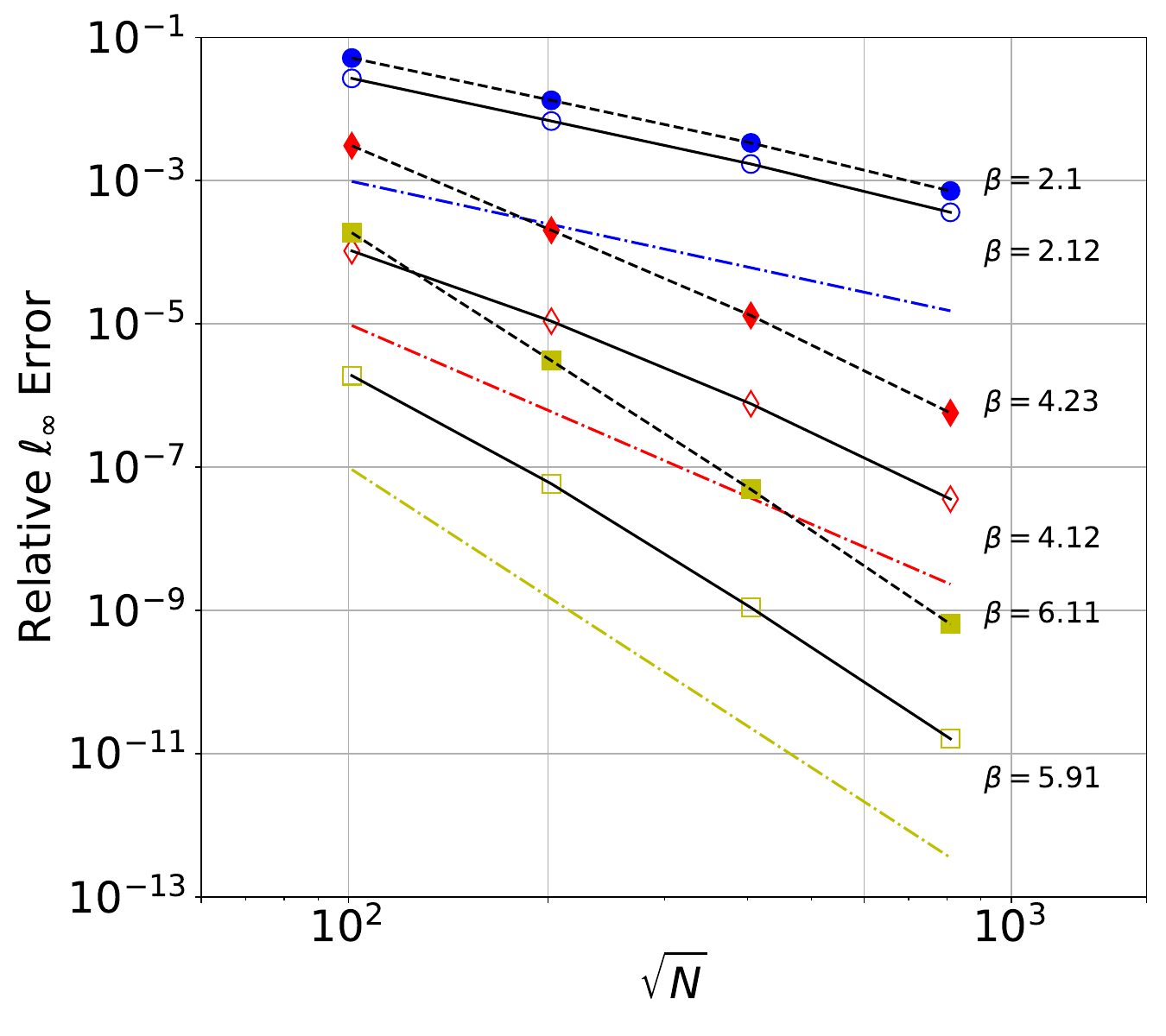} \\ 
\rotatebox{90}{\hspace{0.1\textwidth} \small \quad \textbf{(c) Laplacian}} & \includegraphics[width=0.4\textwidth]{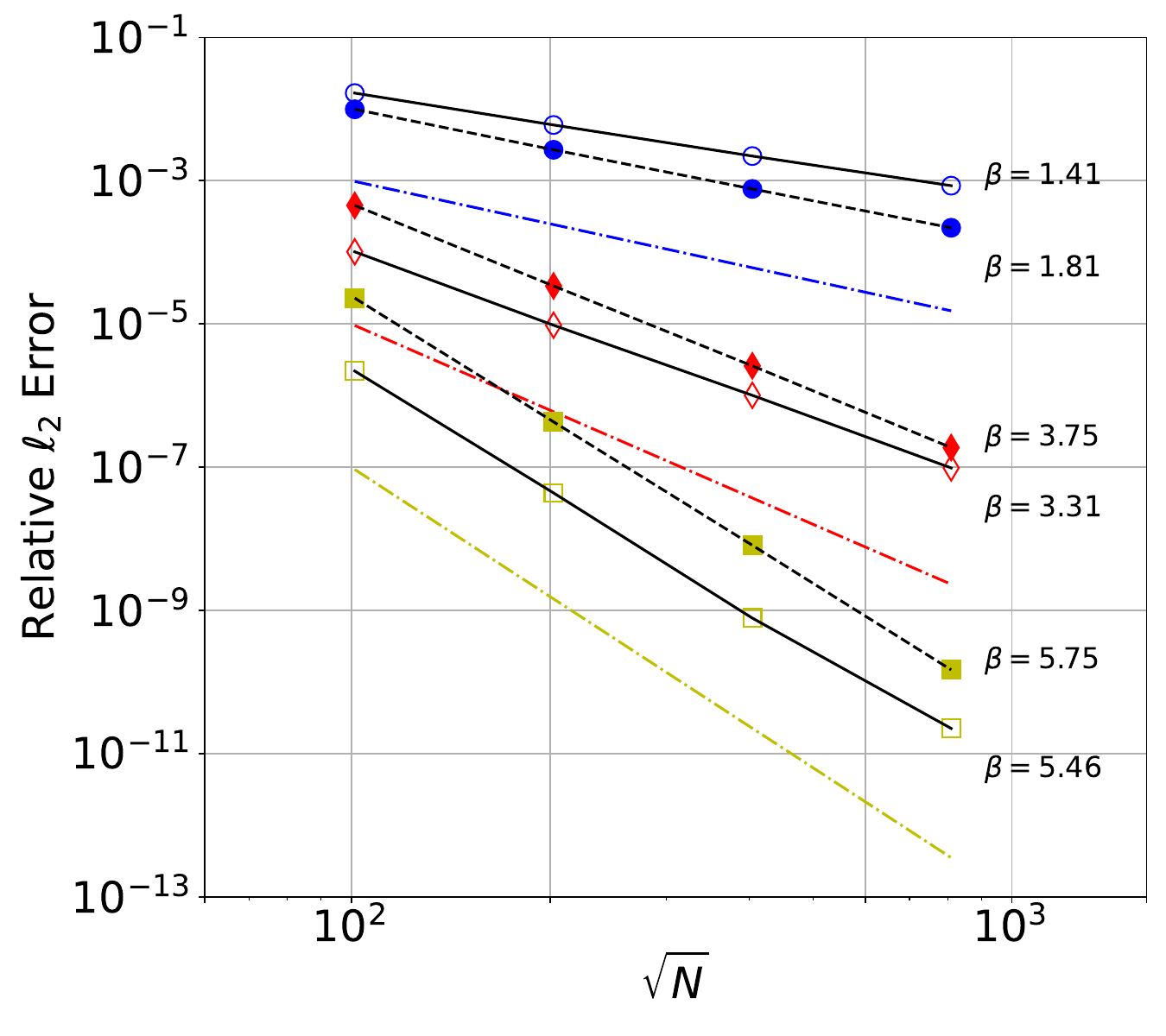} & 
\includegraphics[width=0.4\textwidth]{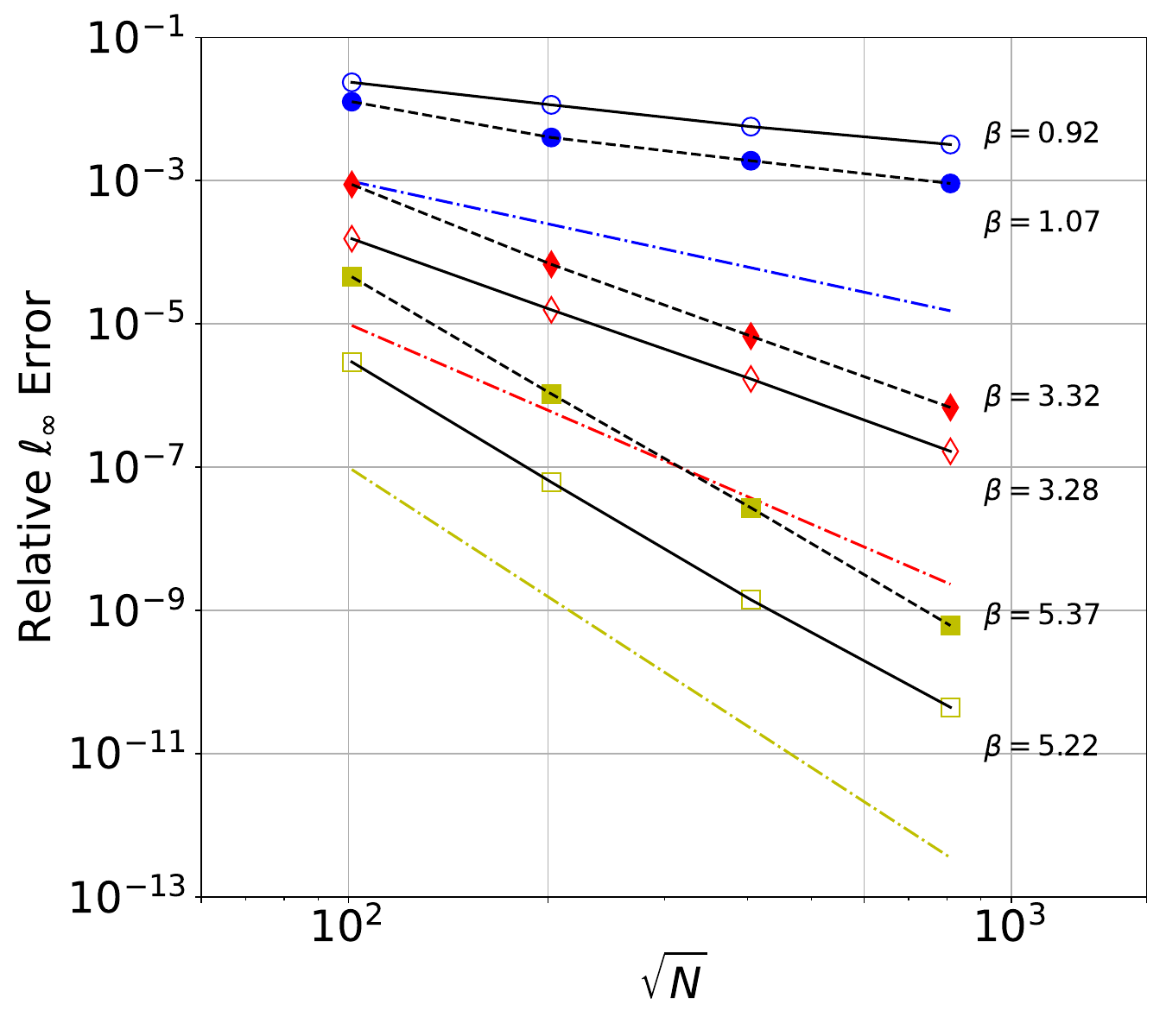} 
\end{tabular}
\caption{Convergence results for (a) surface gradient, (b) divergence, and (b) Laplacian on the sphere using icosahedral node sets.  Errors are given in relative two-norms (first column) and max-norms (second column).  Markers correspond to different $\ell$: filled markers are GMLS and open markers are RBF-FD.  Dash-dotted lines without markers correspond to 2nd, 4th, and 6th order convergence with $1/\sqrt{N}$.  $\beta$ are the measured order of accuracy computed using the lines of best fit to the last three reported errors.  \label{fig:convergence_icos}}
\end{figure}

\begin{figure}
\centering
\begin{tabular}{ccc}
\rotatebox{90}{\hspace{0.1\textwidth}\small \quad \textbf{(a) Gradient}} & \includegraphics[width=0.4\textwidth]{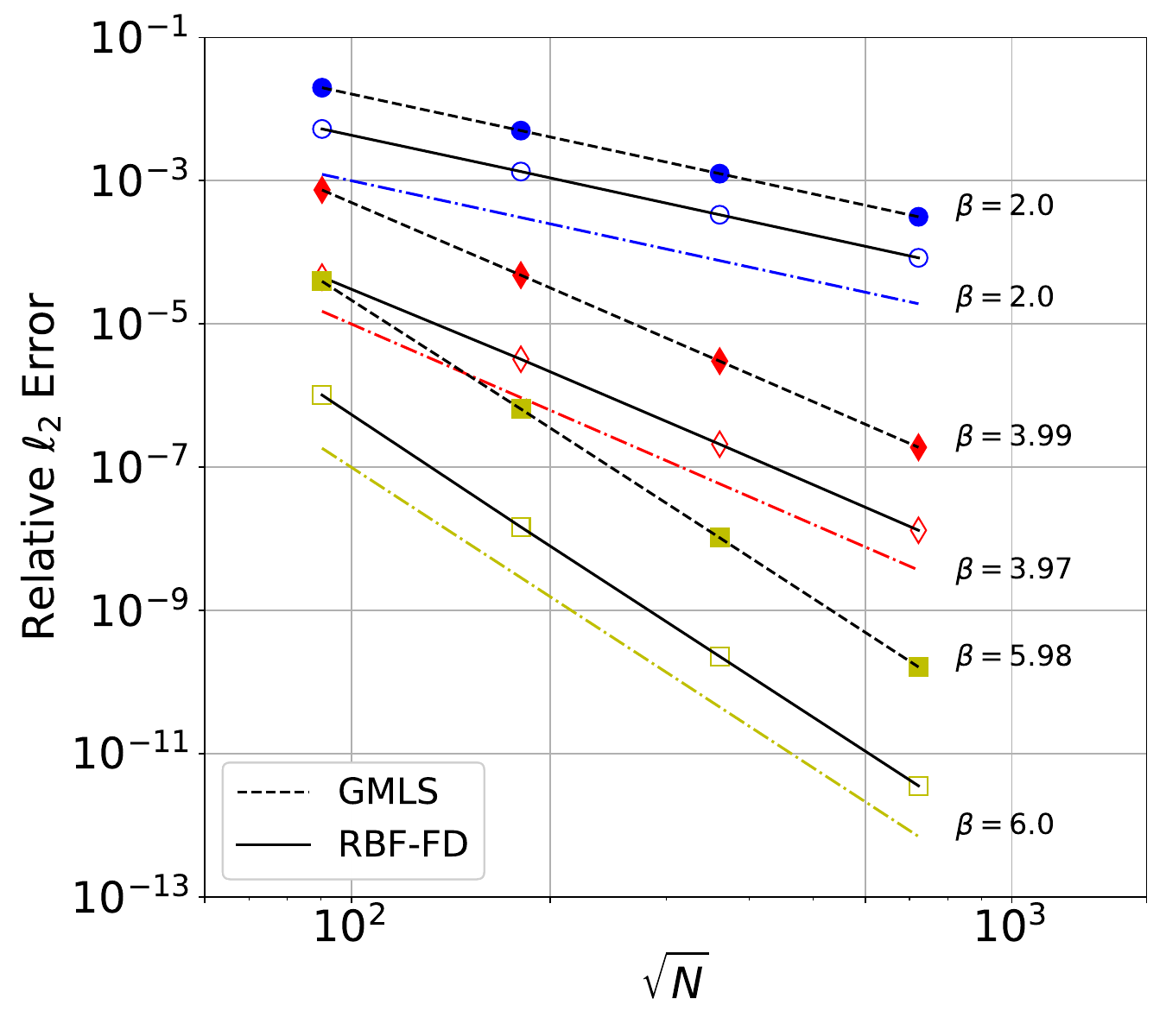} &
\includegraphics[width=0.4\textwidth]{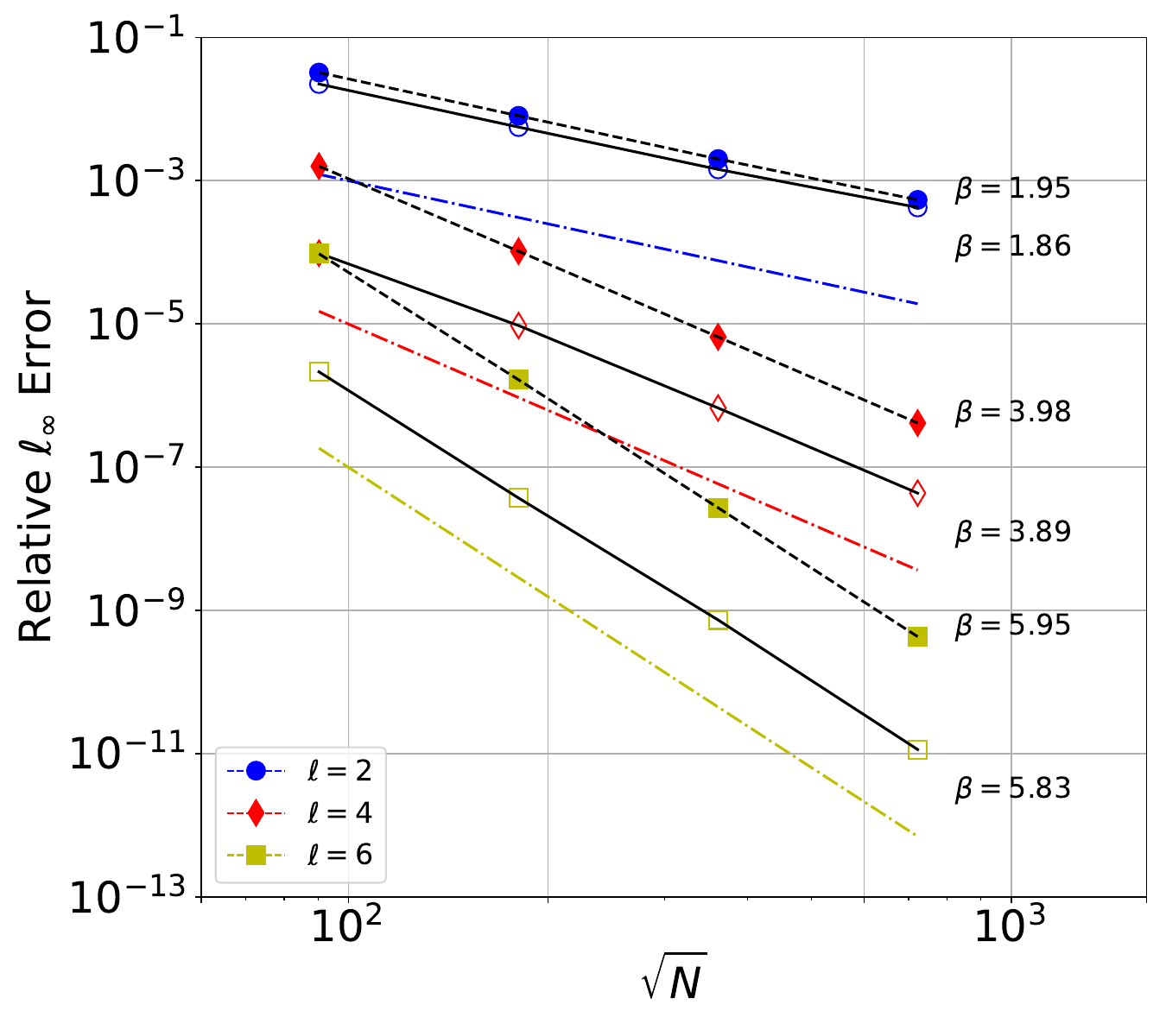} \\
\rotatebox{90}{\hspace{0.1\textwidth} \small \quad \textbf{(b) Divergence}} & \includegraphics[width=0.4\textwidth]{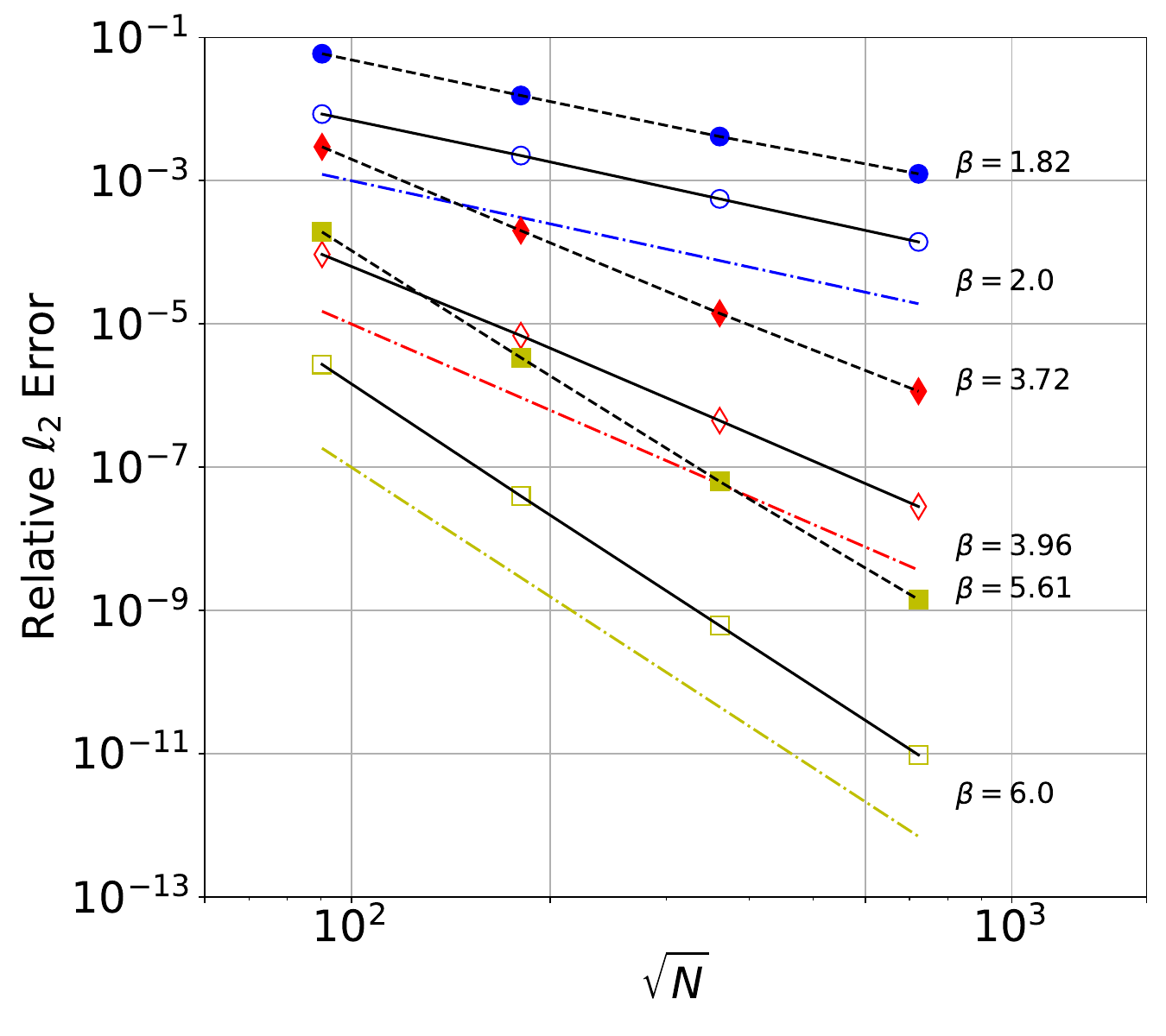} & 
\includegraphics[width=0.4\textwidth]{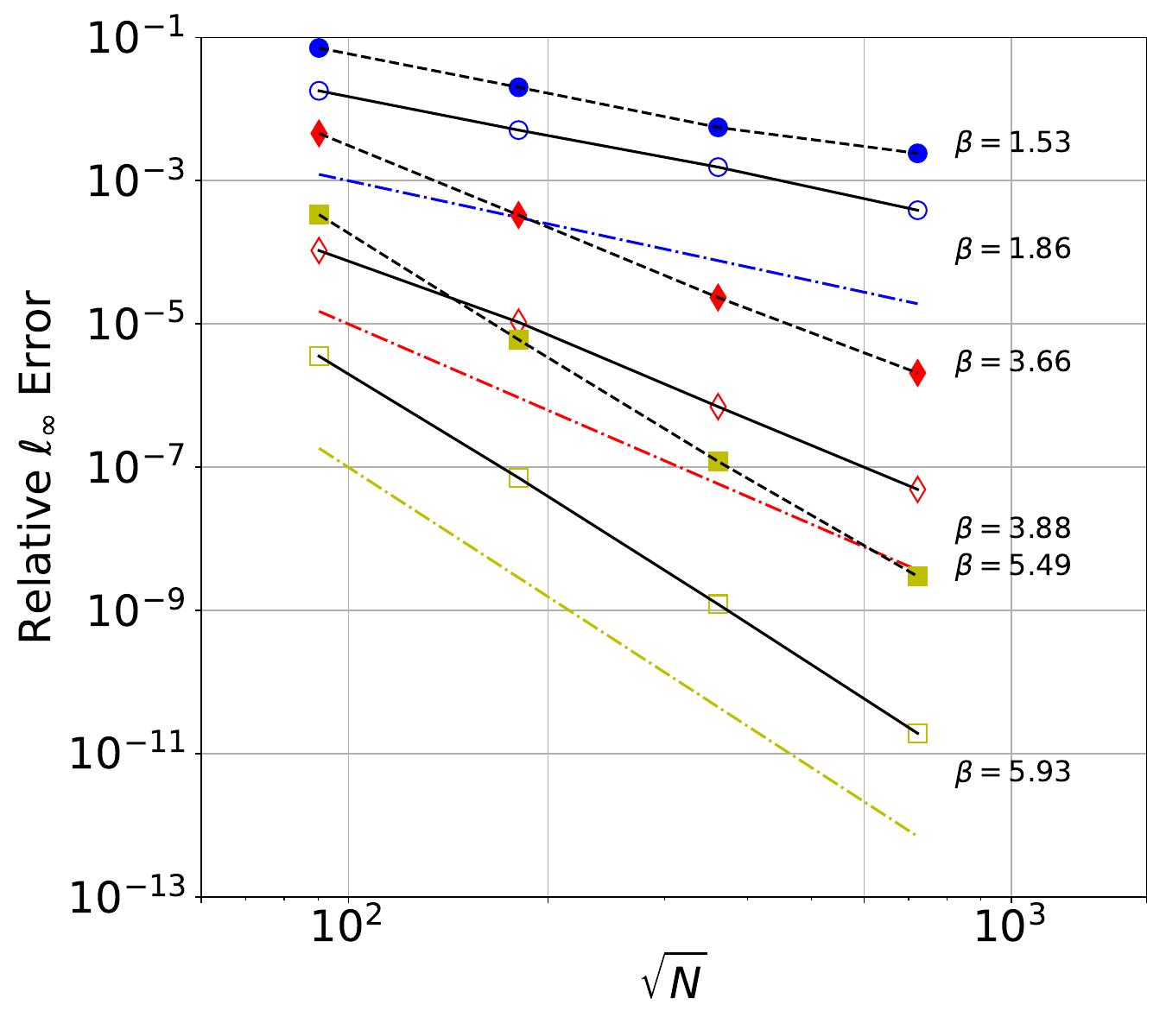} \\ 
\rotatebox{90}{\hspace{0.1\textwidth} \small \quad \textbf{(c) Laplacian}} & \includegraphics[width=0.4\textwidth]{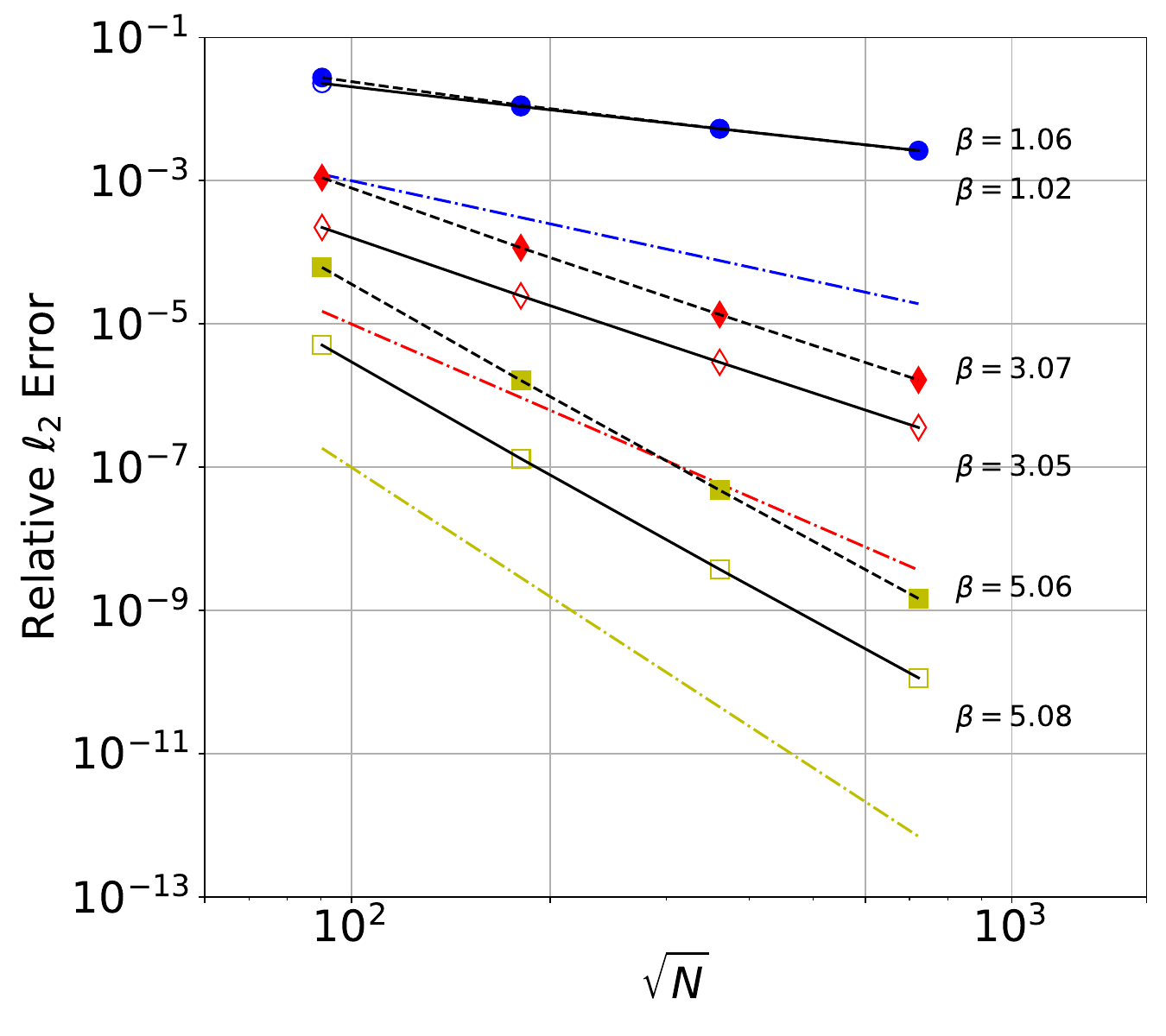} & 
\includegraphics[width=0.4\textwidth]{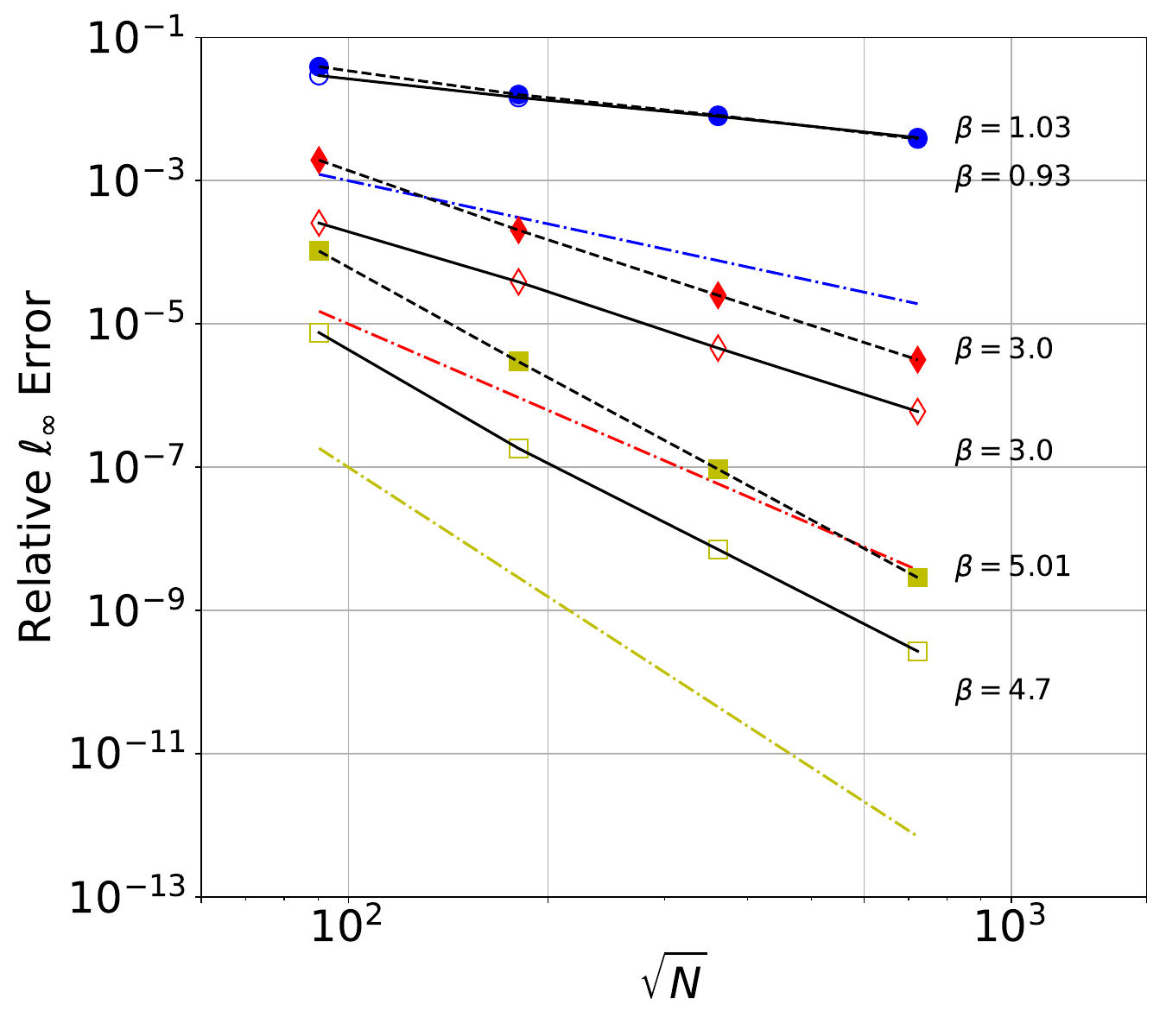} 
\end{tabular}
\caption{\refa{Same as Figure \ref{fig:convergence_icos}, but for the Hammersley nodes on the sphere.}\label{fig:convergence_hammer}}
\end{figure}

Figures \ref{fig:convergence_icos} and \ref{fig:convergence_hammer} display the convergence results for GMLS and RBF-FD as a function of $N$.  Each figure is for a different point set type and contains the results for approximating the surface gradient, divergence, and Laplacian in both the relative two- and max-norms and for different polynomial degrees $\ell$.  We see from all the results that the measured convergence rates for GMLS and RBF-FD are similar, but that RBF-FD gives lower errors for the same $N$ and $\ell$ for approximating the surface gradient and divergence.  \refa{This is also true for the surface Laplacian when $\ell=4$ and $\ell=6$, but not for $\ell=2$.  For this case, GMLS gives lower errors for the same $N$ on the icosahedral nodes and about the same error for the Hammersley nodes.  We also see from Figure \ref{fig:convergence_hammer} that both methods do not appear to be effected by stability issues associated with badly placed points in the stencils for the Hammersley nodes.}  

The measured convergence rates in the two-norm for the surface gradient and divergence approximations are close to the expected rates of $\ell$ for both point sets.  However, when looking at the convergence rates of the surface Laplacian, we see from Figure \ref{fig:convergence_icos} that the icosahedral nodes have higher rates than for the \refa{Hammersley nodes in Figure \ref{fig:convergence_hammer}}.  These improved convergence rates have been referred to as superconvergence in the GMLS literature and rely on the point set being structured so that the stencils have certain symmetries~\cite{LiangZhao13}.  When these symmetries do not exist, as is the case for the \refa{Hammersley nodes}, the convergence rates for the surface Laplacian \refa{more closely follow} the expected rates of $\ell-1$.

\subsection{Convergence comparison: Torus}\label{subsec:torus}
\begin{figure}
\centering
\begin{tabular}{ccc}
\rotatebox{90}{\hspace{0.1\textwidth}\small \quad \textbf{(a) Gradient}} & \includegraphics[width=0.4\textwidth]{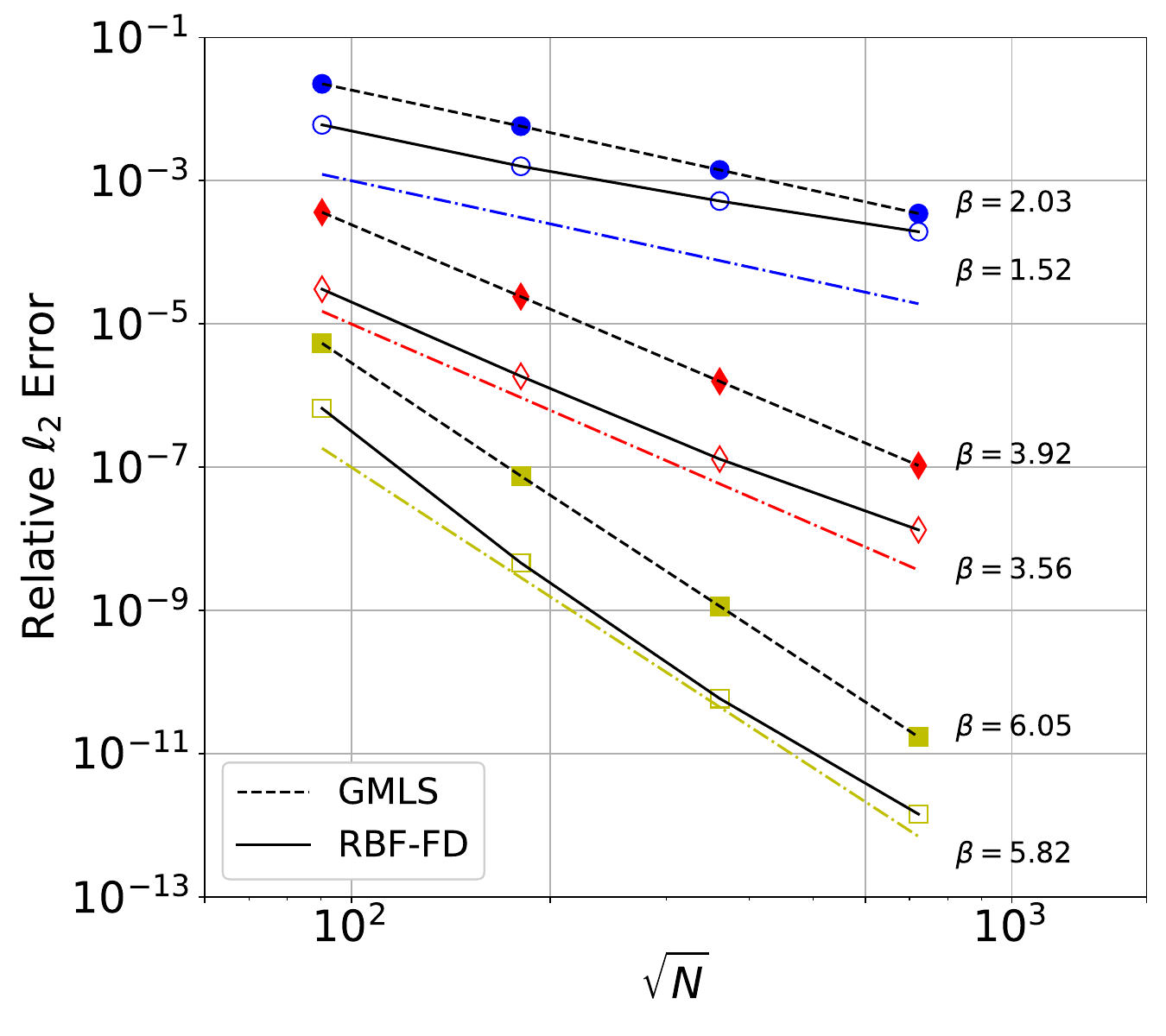} &
\includegraphics[width=0.4\textwidth]{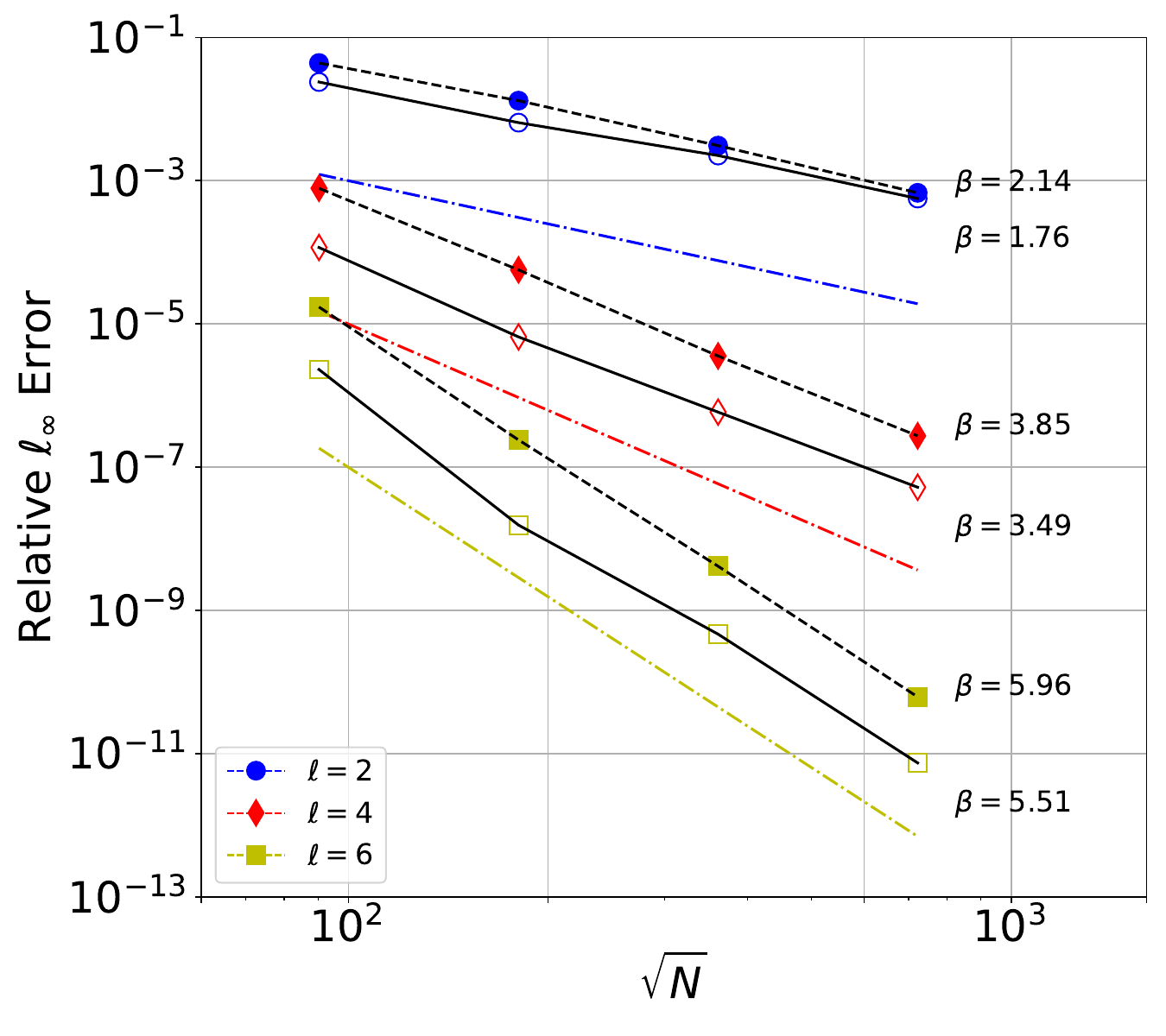} \\
\rotatebox{90}{\hspace{0.1\textwidth} \small \quad \textbf{(b) Divergence}} & \includegraphics[width=0.4\textwidth]{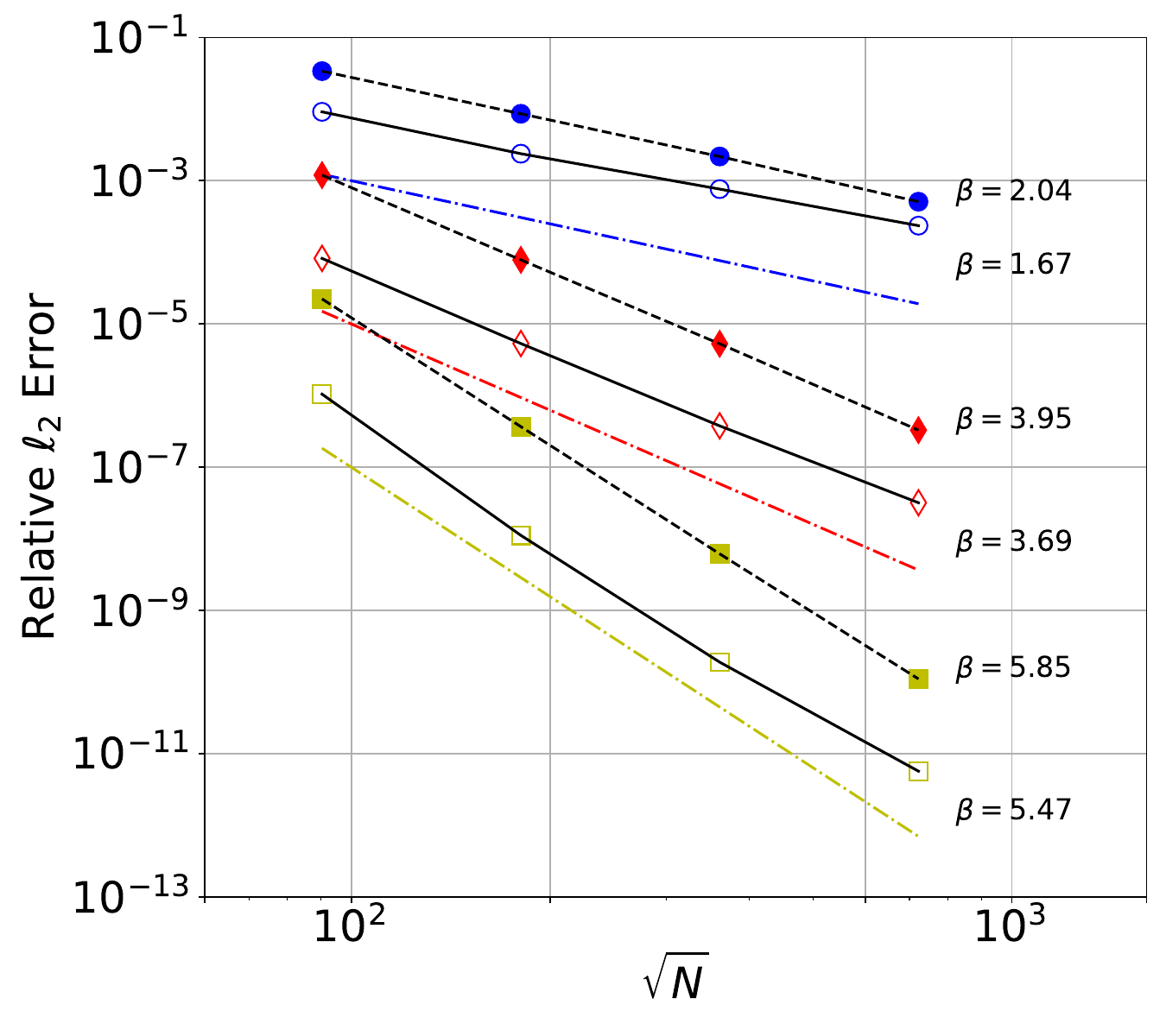} & 
\includegraphics[width=0.4\textwidth]{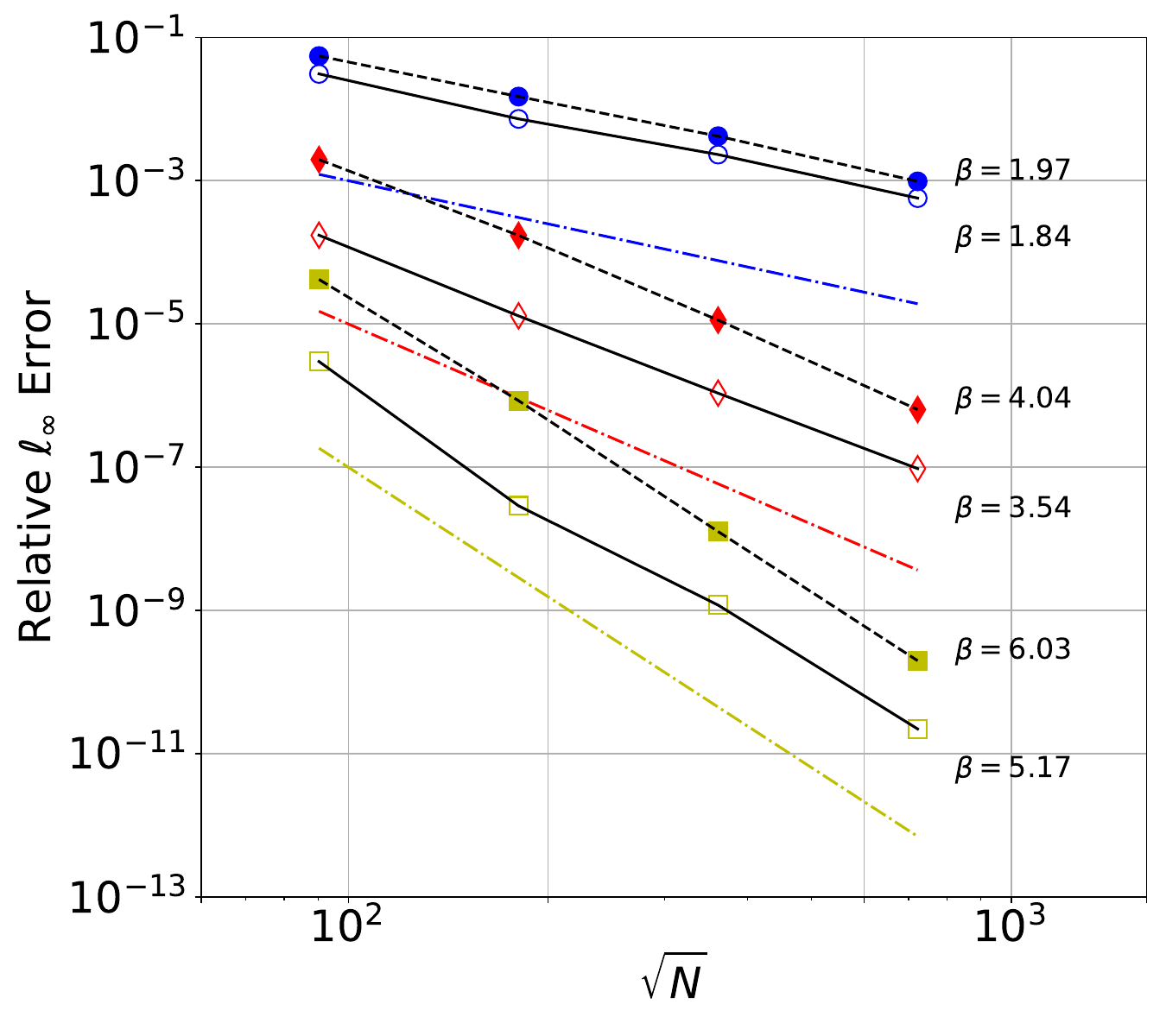} \\ 
\rotatebox{90}{\hspace{0.1\textwidth} \small \quad \textbf{(c) Laplacian}} & \includegraphics[width=0.4\textwidth]{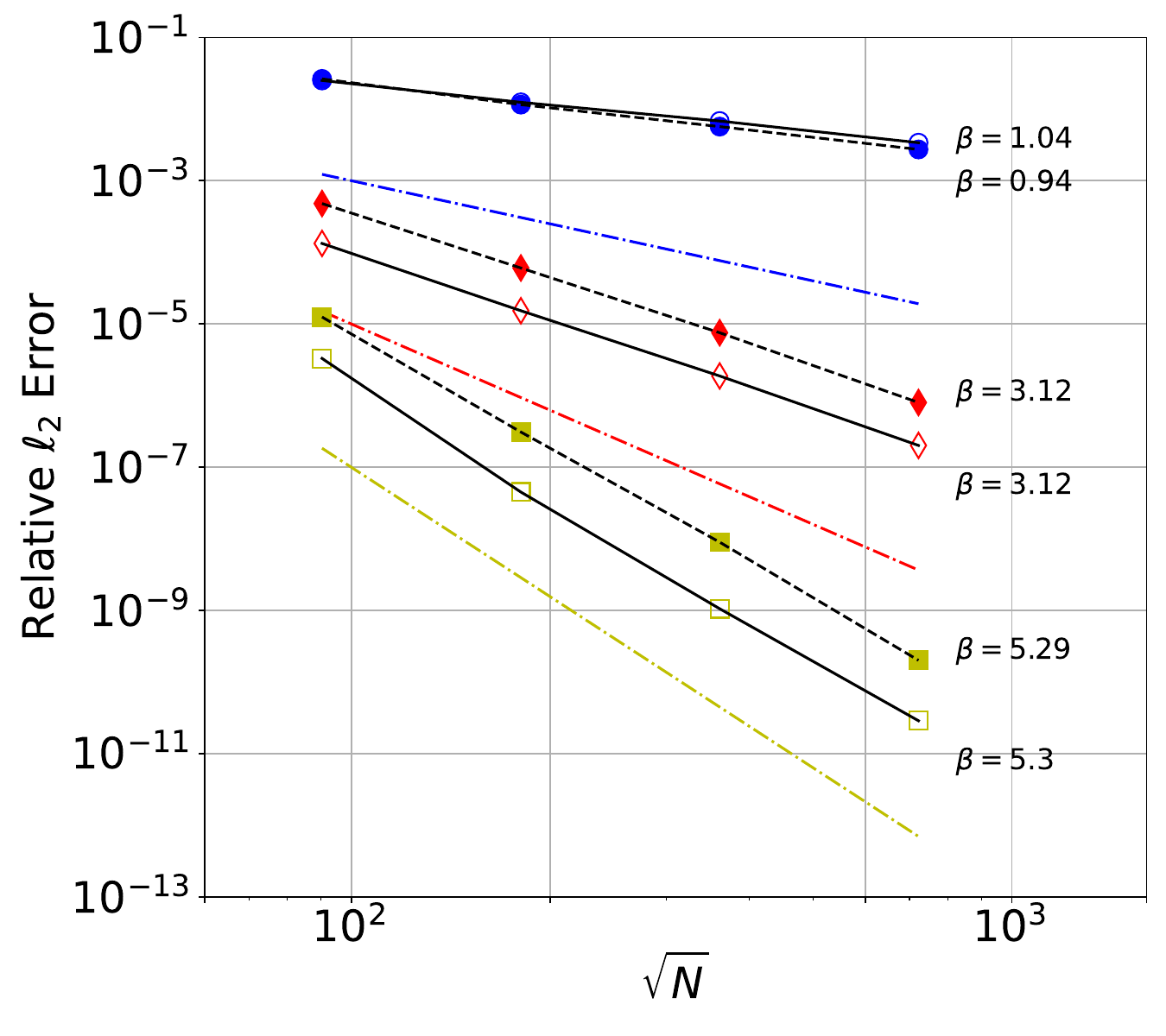} & 
\includegraphics[width=0.4\textwidth]{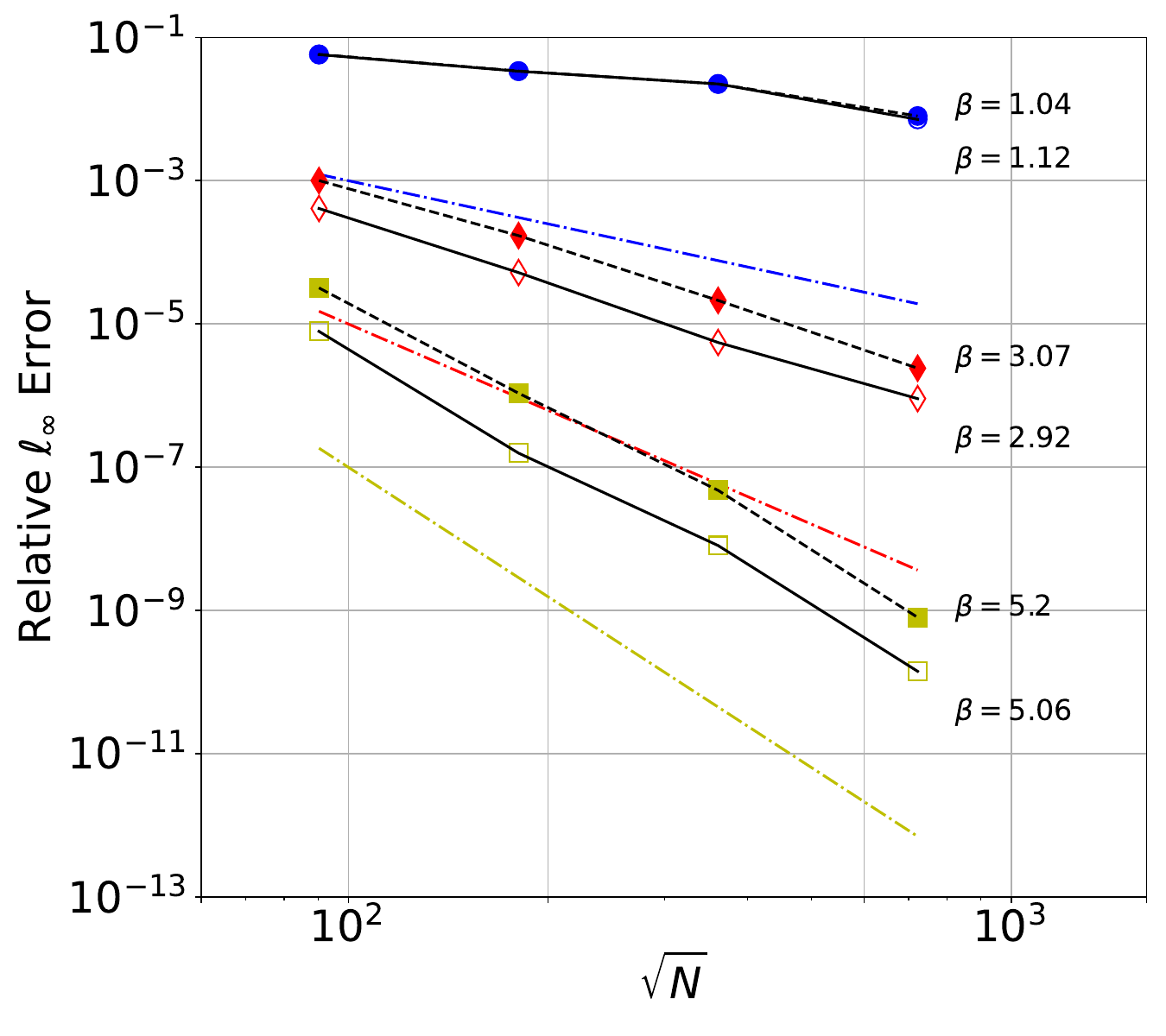} 
\end{tabular}
\caption{Same as Figure \ref{fig:convergence_icos}, but for torus using Poisson disk points.\label{fig:convergence_torus_poisson}}
\end{figure}
The convergence comparisons on the torus are based on the target function
\begin{equation}\label{eq:torus_forcing_exact}
u(\vx) = \ds \frac{x}{8}(x^4 - 10x^2 y^2 + 5y^4)(r^2 - 60 z^2),\; \vx \in \mathbb{T}^2,
\end{equation}
where $r = \sqrt{x^2+y^2}$.  This function has also been used in other studies of RBF methods for surfaces~\cite{FuselierWright2013}.   As with the sphere example, the surface gradient of $u$ can be computed as $\nabla_\M u = \nabla u - \veta(\veta\cdot \nabla u)$, where $\veta$ is the unit outward normal to $\mathbb{T}^2$, which can be computed from the implicit equation \eqref{eq:torus}.  The surface Laplacian of \eqref{eq:torus_forcing_exact} is given in~\cite{FuselierWright2013} as
\begin{align*}
\laps u(\vx) = \ds -\frac{3x}{8r^2}(x^4 - 10x^2 y^2 + 5y^4)(10248r^4 - 34335r^3 + 41359r^2 - 21320r + 4000),\; \vx \in \mathbb{T}^2.
\end{align*}
Similar to the sphere, we use samples of $\nabla_{\M} u$ in the tests of the divergence and compare the results with $\laps u$ above. 

We first study the convergence rates with the stencil radius scaling $\tau=1.5$ and approximate the tangent space, as we did with the sphere tests.  Figure \ref{fig:convergence_torus_poisson} displays the results for the surface gradient, divergence, and Laplacian.  We see that errors for RBF-FD are again smaller than the errors for GMLS in almost all cases over the range of $N$ tested.  However, GMLS has a slightly higher convergence rates in the case of the surface gradient and divergence, but not for the Laplacian.  Both methods have convergence rates that are close to the expected rates of $\ell$ for these surface gradient and divergence and $\ell-1$ for the Laplacian.

\begin{figure}
\centering
\begin{tabular}{cc}
\includegraphics[width=0.4\textwidth]{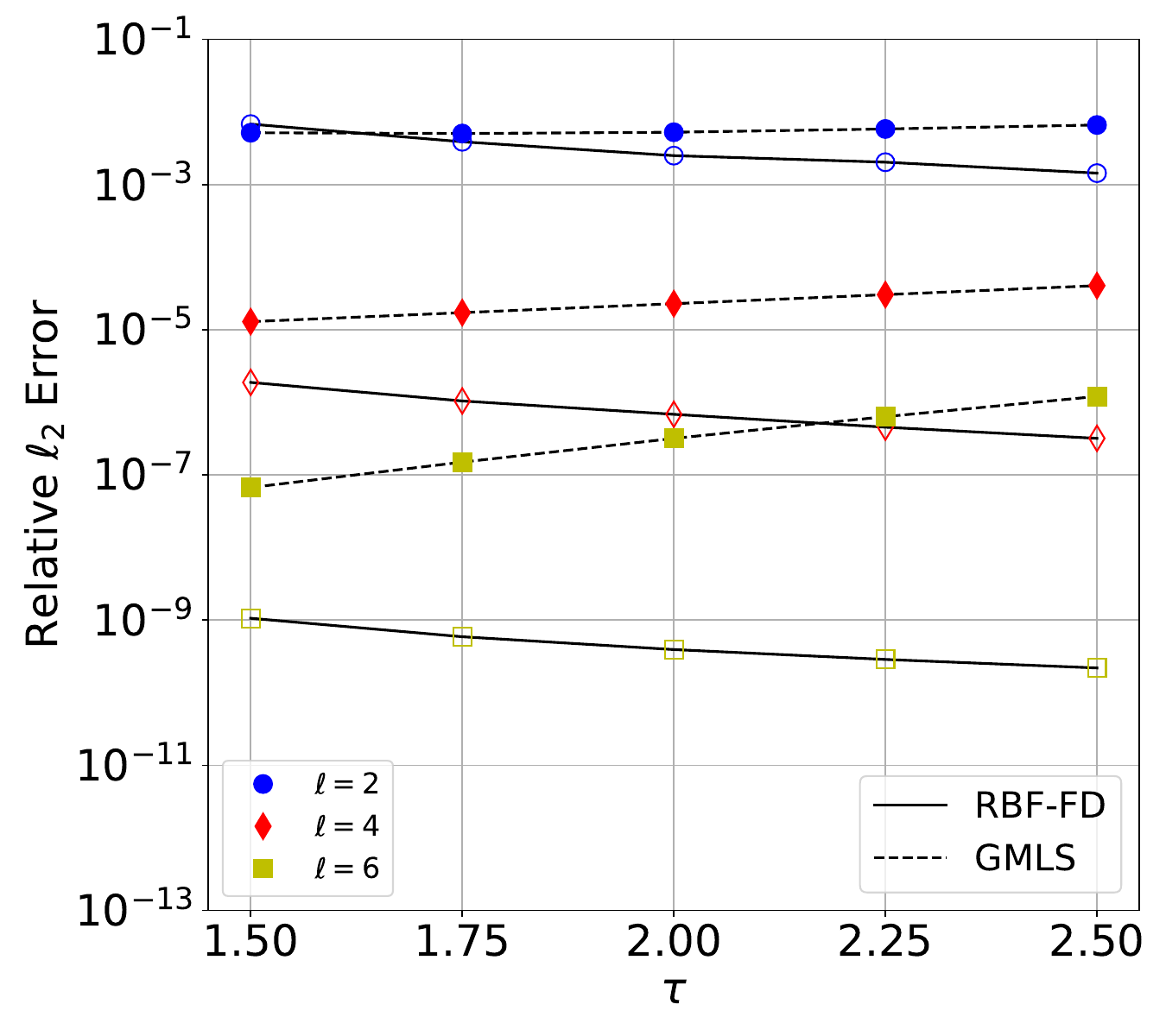} &
\includegraphics[width=0.4\textwidth]{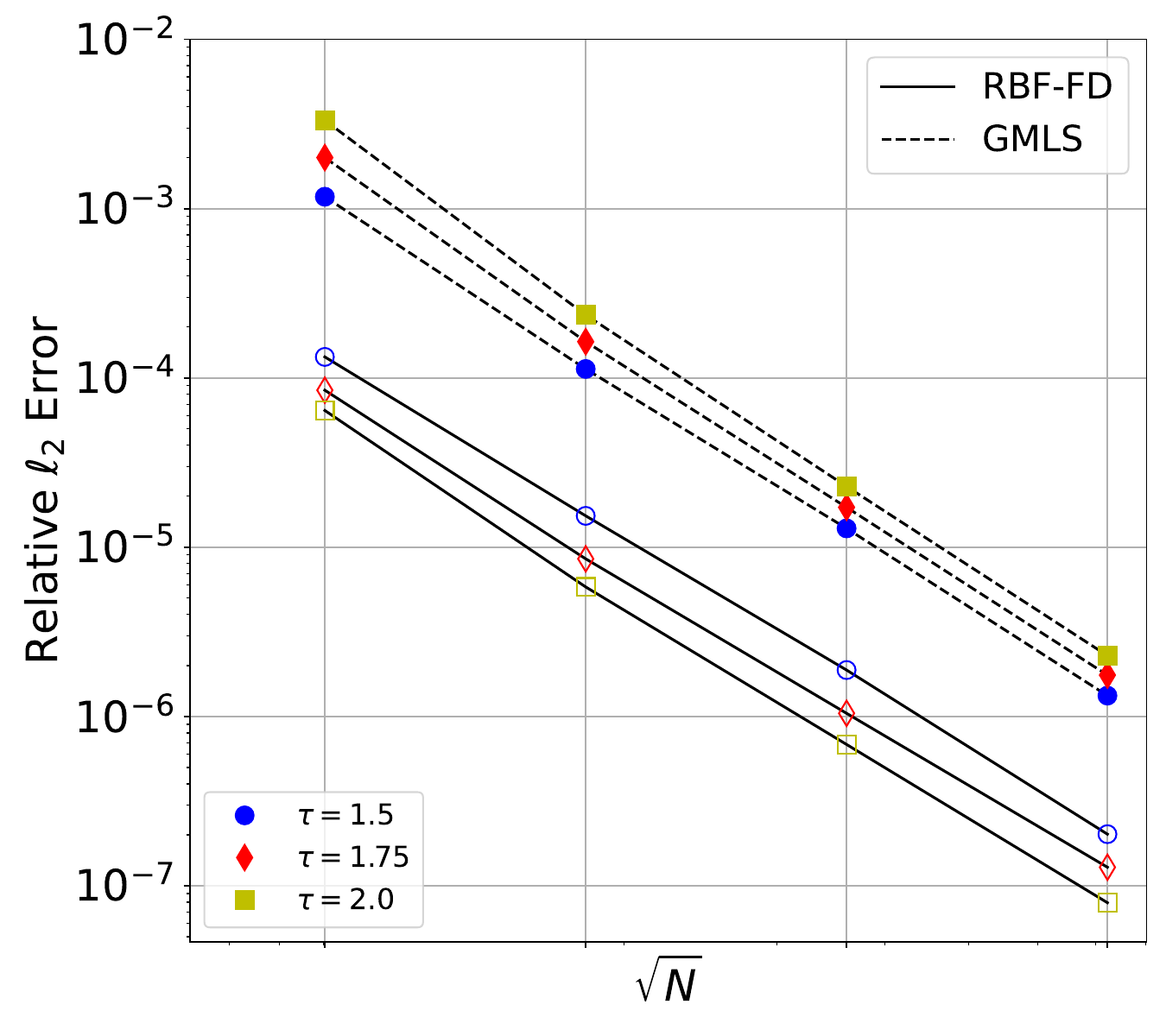}
\end{tabular}
\caption{Relative two-norm errors of the surface Laplacian approximations as the stencil radius parameter $\tau$ varies.  Left figure shows errors for several different values of $\tau$ and a fixed $N=130463$.  Right figure shows the convergence rates of the methods for different $\tau$ and a fixed $\ell=4$.\label{fig:tau_scaling}}
\end{figure}

Next we investigate how the approximation properties of the two methods change when $\tau$ is increased, which results in larger stencil sizes. We focus on approximating the surface Laplacian as similar results were found for the other SDOs.  In the left plot of Figure \ref{fig:tau_scaling}, we show the relative two-norm errors of the approximations for a fixed $N$ as $\tau$ varies from 1.5 to 2.5.  We see that increasing $\tau$ has opposite effects on the two methods: the errors decrease for RBF-FD and increase with GMLS.  We see similar results in the right plot of Figure \ref{fig:tau_scaling}, where we show the convergence of the methods with increasing $N$ for different fixed values of $\tau$ (and $\ell$ fixed at $4$).  While the convergence rates do not appear to change with $\tau$, the overall errors decrease for RBF-FD and increase for GMLS.  \refa{It should be noted that the errors eventually increase for GMLS as $\tau$ decreases to 1 (which has been observed in other studies) and picking an optimal $\tau$ in an automated way is challenging (e.g.~\cite{lipman2006error,wang2008optimal}).}

These results make sense when one considers the different types of approximations the methods are based on: RBF-FD is based on interpolation, while GMLS is based on least squares approximation.  As the stencil sizes increase, RBF-FD has a larger approximation space consisting of more shifts of PHS kernels, which can reduce the errors~\cite{davydov2019optimal}.  However, GMLS has the same fixed approximation space of polynomials of degree $\ell$ regardless of the stencil size.  

Finally, we compare the errors when the exact and approximate tangent spaces are used in the two methods.  We focus only on the surface Laplacian and for $\ell=4$ since similar results were obtained for the other operators and other $\ell$.  Table \ref{tbl:tangent_space} shows the results for both methods.  The approximate tangent spaces were computed using the methods from Sections \ref{sec:gmls_tangent_space} (GMLS) and \ref{sec:rbffd_tangent_space} (RBF-FD) also using the polynomial degree $\ell=4$.  \refa{As discussed in Section 5, this choice is made so that the tangent spaces are approximated with the same asymptotic order of accuracy as the approximation of the metric terms with GMLS.}  We see from the table that the differences between using the exact or the approximate tangent spaces for approximating the surface Laplacian is minor.
\begin{center}
\begin{table}[H]
\centering\begin{tabular}{|c||cc|cc|}
\hline
 & \multicolumn{2}{c|}{GMLS} & \multicolumn{2}{c|}{RBF-FD} \\
\hline
$N$ & Exact & Approx. & Exact & Approx. \\ 
\hline
\hline
 8153 & 4.7984e-04 & 4.8004e-04 & 1.3311e-04 & 1.3312e-04 \\
\hline
 32615 & 6.0457e-05 & 6.04654e-05 & 1.5321e-05 & 1.5322e-05 \\
\hline
 130463 & 7.5486e-06 & 7.5488e-06 & 1.8811e-06 & 1.8811e-06 \\
\hline
 521855 & 8.0158e-07 & 8.0159e-07 & 2.0177e-07 & 2.0176e-07 \\
\hline
\end{tabular}
\caption{Comparison of the relative $\ell_2$ errors for the surface Laplacian on the torus using the exact tangent space for the torus and approximations to it based on the methods from Sections \ref{sec:gmls_tangent_space} (GMLS) and \ref{sec:rbffd_tangent_space} (RBF-FD). In all cases, $\ell=4$ and the points are based on Poisson disk sampling.\label{tbl:tangent_space}}
\end{table}
\end{center}

\subsection{Efficiency comparison}
\begin{figure}
\centering
\begin{tabular}{cc}
\centering
\includegraphics[width=0.45\textwidth]{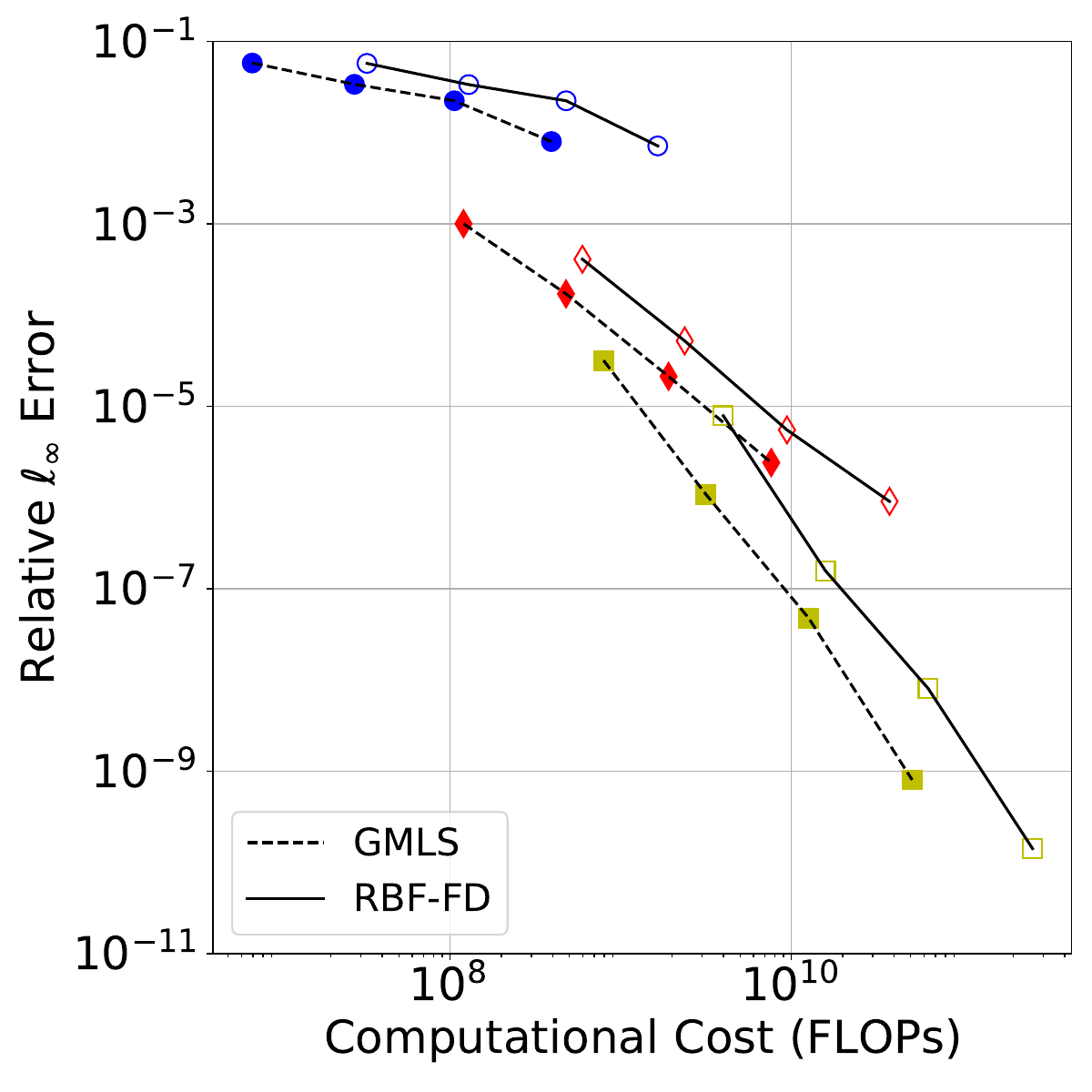} & \includegraphics[width=0.45\textwidth]{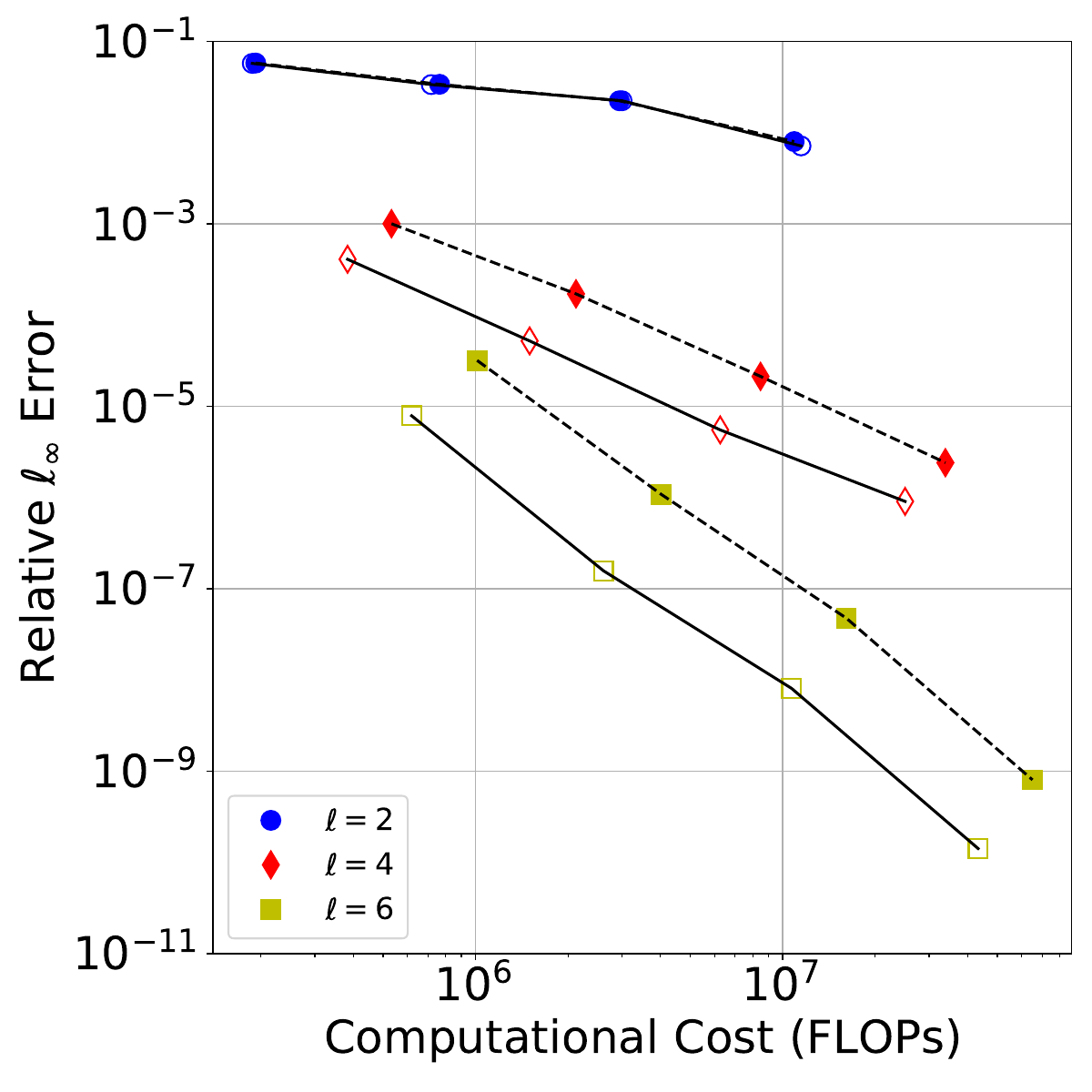} \\
(a) \refa{Set-up + Evaluation} & (b) Evaluation
\end{tabular}
\caption{\refa{Comparison of the computational efficiency of GMLS and RBF-FD for approximating the surface Laplacian in terms of accuracy per computational cost.  (a) Shows the efficiency when considering the setup and evaluation costs as defined in \eqref{eq:setup_cost} and \eqref{eq:runtime_cost}, respectively, while (b) show the efficiency when only considering the evaluation cost. \label{fig:buildcost_figs}}}
\end{figure}

The results in Section \ref{sec:results} demonstrate that RBF-FD and GMLS have similar asymptotic convergence rates for the same $\ell$, but that RBF-FD can achieve lower errors for the same $N$ and stencil sizes.  In this section, we consider which of the \refa{methods are} more computationally efficient in terms of error per computational cost.  We examine both the efficiency when the setup costs are included and when just the evaluation costs are included, as measured by \eqref{eq:setup_cost} and \eqref{eq:runtime_cost}, respectively.  We limit this comparison to $\tau=1.5$, but note that it may be possible to tune this parameter to (marginally) optimize the efficiency of either method over this case.  Figure \ref{fig:buildcost_figs} displays the results of this examination for the case of computing the surface Laplacian on the \refc{torus} discretized with Poisson disk sampling.  Similar results were obtained for other SDOs and for the sphere, so we omit them.  We see from the figure that GMLS is more efficient when the setup costs are included, but that RBF-FD is more efficient when only evaluation costs are included.  For problems where the point sets are fixed and approximations to a SDO are required to be performed multiple times---as occurs when solving a time-dependent surface PDEs---the setup costs are not as important as the evaluation costs since they are amortized across all time-steps.  In this scenario RBF-FD is the more efficient method.

\section{Concluding remarks}\label{sec:remarks}
We presented a thorough comparison of the GMLS and RBF-FD methods for approximating the three most common SDOs: the gradient, divergence, and Laplacian (Laplace-Beltrami).  Our analysis of the two different formulations of SDOs used in the methods revealed that if the exact tangent space for the surface is used, these formulations are identical.  We further derived a new RBF-FD method for approximating the tangent space of surfaces represented only by point clouds.  Our numerical investigation of the methods showed that they appear to converge at similar rates when the same polynomial degree $\ell$ is used, but that RBF-FD generally gives lower errors for the same $N$ and $\ell$.  We additionally examined the dependency of the stencil size on the methods (as measured by the $\tau$ parameter) and found that the errors produced by GMLS deteriorate as the stencil size increases.  The errors for RBF-FD, contrastingly, appear to keep improving as the stencil size increases.  However, we don't expect this trend to continue indefinitely, as eventually the tangent plane formulation breaks down when the stencil size becomes too large. Finally, we investigated the computational efficiency of the methods in terms of error versus computational cost and found GMLS to be more efficient when setup costs are included and RBF-FD to be more efficient when only considering evaluation costs. 



\section*{Acknowledgements} 
\paragraph{Funding}
AMJ was partially supported by US NSF grant CCF-1717556.  PAB \& PAK were supported by the U.S.\ Department of Energy, Office of Science, Advanced Scientific Computing Research (ASCR) Program and Biological and Environmental Research (BER) Program under a Scientific Discovery through Advanced Computing (SciDAC 4) BER partnership pilot project.   PAK was additionally supported by the Laboratory Directed Research \& Development (LDRD) program at Sandia National Laboratories and ASCR under Award Number DE-SC-0000230927. AMJ was also partially supported by the Climate Model Development and Validation (CMDV) program, funded by BER.  Part of this work was conducted while AMJ was employed at the Computer Science Research Institute at Sandia National Laboratories.  GBW was partially supported by U.S.\ NSF grants CCF-1717556 and DMS-1952674.



\bibliographystyle{elsarticle-num} 
\bibliography{biber}

\section*{Statements and Declarations}
Sandia National Laboratories is a multi-mission laboratory managed and operated by National Technology \& Engineering Solutions of Sandia, LLC (NTESS), a wholly owned subsidiary of Honeywell International Inc., for the U.S. Department of Energy’s National Nuclear Security Administration (DOE/NNSA) under contract DE-NA0003525. This written work is authored by an employee of NTESS. The employee, not NTESS, owns the right, title and interest in and to the written work and is responsible for its contents. Any subjective views or opinions that might be expressed in the written work do not necessarily represent the views of the U.S. Government. The publisher acknowledges that the U.S. Government retains a non-exclusive, paid-up, irrevocable, world-wide license to publish or reproduce the published form of this written work or allow others to do so, for U.S. Government purposes. The DOE will provide public access to results of federally sponsored research in accordance with the DOE Public Access Plan.

%

\paragraph{Competing Interests}

The authors have no relevant financial or non-financial interests to disclose.

\paragraph{Author Contributions}

\begin{itemize}
\item AMJ: Conceptualization, Methodology, Software, Formal analysis, Investigation, Writing: Original draft preparation, Writing: Review \& Editing.

\item PAB: Methodology, Resources, Writing: Review \& Editing, Funding acquisition, Supervision

\item PAK: Resources, Software, Writing: Review \& Editing

\item GBW: Methodology, Writing: Orignal draft preparation, Writing: Review \& Editing, Funding acquisition, Supervision

\end{itemize}

%
\end{document}